\documentclass{amsart}

\usepackage{amssymb} \usepackage{amsthm} \usepackage{amsmath} \usepackage{comment,color}
\usepackage[hidelinks]{hyperref}

\newcommand{\numberseries}{\bfseries}   
\newlength{\thmtopspace}                
\newlength{\thmbotspace}                
\newlength{\thmheadspace}               
\newlength{\thmindent}                  
\newtheoremstyle{bfupright head,slanted body}
                {\thmtopspace}{\thmbotspace}
                {\slshape}{\thmindent}{\bfseries}{.}{\thmheadspace}
                {{\numberseries \thmnumber{(#2) }}\thmnote{#3}}
\newtheoremstyle{fixed bf head,slanted body}
                {\thmtopspace}{\thmbotspace}{\slshape}
                {\thmindent}{\bfseries}{.}{\thmheadspace}
                {{\numberseries \thmnumber{#2\;}}\thmname{#1}\thmnote{ (#3)}}
\newtheoremstyle{fixed bf head,upright body}
                {\thmtopspace}{\thmbotspace}{\upshape}
                {\thmindent}{\bfseries}{.}{\thmheadspace}
                {{\numberseries \thmnumber{#2\;}}\thmname{#1}\thmnote{ (#3)}}
\newtheoremstyle{numbered paragraph}
                {\thmtopspace}{\thmbotspace}{\upshape}
                {\thmindent}{\upshape}{}{\thmheadspace}
                {{\numberseries \thmnumber{#2.}}}
\theoremstyle{bfupright head,slanted body}
\newtheorem{res}{}[section]             \newtheorem*{res*}{}
\theoremstyle{fixed bf head,slanted body}
\newtheorem{thm}[res]{Theorem}          \newtheorem*{thm*}{Theorem}
\newtheorem{prp}[res]{Proposition}      \newtheorem*{prp*}{Proposition}
\newtheorem{cor}[res]{Corollary}        \newtheorem*{cor*}{Corollary}
            \newtheorem*{lem*}{Lemma}
\theoremstyle{fixed bf head,upright body}
\newtheorem{exa}[res]{Example}          \newtheorem*{exa*}{Example}
\newtheorem{rmk}[res]{Remark}
\theoremstyle{numbered paragraph}
\newtheorem{ipg}[res]{}
\newlength{\thmlistleft}        
\newlength{\thmlistright}       
\newlength{\thmlistpartopsep}   
\newlength{\thmlisttopsep}      
\newlength{\thmlistparsep}      
\newlength{\thmlistitemsep}     
\newcommand{\eqclbl}[1]{{\upshape(\textit{#1})}}
\newcounter{prt}
\newenvironment{prt}{\begin{list}{\upshape (\alph{prt})}%
    {\usecounter{prt}%
      \setlength{\leftmargin}{\thmlistleft}%
      \setlength{\labelwidth}{\thmlistleft}%
      \setlength{\rightmargin}{\thmlistright}%
      \setlength{\partopsep}{\thmlistpartopsep}%
      \setlength{\topsep}{\thmlisttopsep}%
      \setlength{\parsep}{\thmlistparsep}%
      \setlength{\itemsep}{\thmlistitemsep}}}%
  {\end{list}}%
\newenvironment{prf*}[1][Proof]{%
  \begin{proof}[\bf #1]
    \setcounter{equation}{0}
    }
  {\end{proof}
}
\newcommand{\proofoftag}[2][:]{(#2)#1}\newcommand{\pgref}[1]{\ref{#1}}
\newcommand{\thmref}[2][Theorem~]{#1\pgref{thm:#2}}
\newcommand{\corref}[2][Corollary~]{#1\pgref{cor:#2}}
\newcommand{\prpref}[2][Proposition~]{#1\pgref{prp:#2}}

\newcommand{\exaref}[2][Example~]{#1\pgref{exa:#2}}
\newcommand{\rmkref}[2][Remark~]{#1\pgref{rmk:#2}}
\renewcommand{\eqref}[1]{(\pgref{eq:#1})}
\newcommand{\thmcite}[2][?]{\cite[Theorem~#1]{#2}}
\newcommand{\corcite}[2][?]{\cite[Corollary~#1]{#2}}
\newcommand{\prpcite}[2][?]{\cite[Proposition~#1]{#2}}
\numberwithin{equation}{res}
\def\urltilda{\kern -.15em\lower .7ex\hbox{\~{}}\kern .04em} 
\newcommand{\set}[2][\mspace{1mu}]{\{#1 #2 #1\}}
 \newcommand{\setof}[3][\mspace{1mu}]{\{#1#2 \mid #3#1\}}
\newcommand{\ZZ}{\mathbb{Z}}
\newcommand{\qtext}[1]{\quad\text{#1}\quad}
\newcommand{\qqtext}[1]{\qquad\text{#1}\qquad}
\newcommand{\qand}{\qtext{and}}
\newcommand{\qqand}{\qqtext{and}}
\newcommand{\deq}{\:=\:}
\newcommand{\dge}{\:\ge\:}
\newcommand{\dle}{\:\le\:}
\newcommand{\dqis}{\:\qis\:}
\newcommand{\fa}{\mathfrak{a}}
\newcommand{\fm}{\mathfrak{m}}
\newcommand{\fp}{\mathfrak{p}}
\newcommand{\fq}{\mathfrak{q}}
\newcommand{\is}{\cong}
\newcommand{\qis}{\simeq}
\renewcommand{\le}{\leqslant}
\renewcommand{\ge}{\geqslant}
\newcommand{\Rhat}{\widehat{R}}
\renewcommand{\H}[2][]{\operatorname{H}_{#1}(#2)}
\newcommand{\amp}[1]{\operatorname{amp}#1}
\newcommand{\SpecR}{\operatorname{Spec}R}
\newcommand{\dptR}{\operatorname{depth}R}
\newcommand{\dimR}{\operatorname{dim}R}
\newcommand{\cmdR}{\operatorname{cmd}R}
\newcommand{\Tor}[4][R]{\operatorname{Tor}^{#1}_{#2}(#3,#4)}
\newcommand{\supp}[2][R]{\operatorname{supp}_{#1}#2}
\newcommand{\cosupp}[2][R]{\operatorname{cosupp}_{#1}#2}
\newcommand{\E}[2][R]{\operatorname{E}_{#1}(#2)}
\renewcommand{\dim}[2][R]{\operatorname{dim}_{#1}#2}
\newcommand{\cmd}[2][R]{\operatorname{cmd}_{#1}#2}
\newcommand{\wdt}[2][R]{\operatorname{width}_{#1}#2}
\newcommand{\dpt}[2][R]{\operatorname{depth}_{#1}#2}
\newcommand{\fd}[2][R]{\operatorname{fd}_{#1}#2}
\newcommand{\id}[2][R]{\operatorname{id}_{#1}#2}
\newcommand{\pd}[2][R]{\operatorname{pd}_{#1}#2}
\newcommand{\Hom}[3][R]{\operatorname{Hom}_{#1}(#2,#3)}
\newcommand{\RHom}[3][R]{\operatorname{\mathbf{R}Hom}_{#1}(#2,#3)}
\newcommand{\tp}[3][R]{\nobreak{#2\otimes_{#1}#3}}
\newcommand{\tpp}[3][R]{(\tp[#1]{#2}{#3})}
\newcommand{\Ltp}[3][R]{\nobreak{#2\otimes_{#1}^{\mathbf{L}}#3}}
\newcommand{\Ltpp}[3][R]{(\Ltp[#1]{#2}{#3})}
\newcommand{\RGam}[2][\mathfrak{a}]{\nobreak{\mathbf{R}\Gamma_{#1}#2}}
\newcommand{\lc}[3][\fa]{\operatorname{H}_{#1}^{#2}(#3)}
\newcommand{\LLam}[2][\mathfrak{a}]{\nobreak{\mathbf{L}\Lambda^{#1}#2}}
\newcommand{\catba}{\sqsubset}
\newcommand{\catbb}{\sqsupset}
\newcommand{\catb}{\sqsubset\mspace{-13mu}\sqsupset}
\newcommand{\Cat}[2]{{\mathsf{#2}}(#1)}
\newcommand{\Catsup}[3]{{\mathsf{#2}}^{\text{\upshape #3}}(#1)}
\newcommand{\Catsub}[3]{{\mathsf{#2}}_{#3}(#1)}
\newcommand{\Catsupsub}[4]{{\mathsf{#2}}^{\text{\upshape #3}}_{#4}(#1)}
\newcommand{\D}[1][R]{\Cat{#1}{D}}
\newcommand{\Dba}[1][R]{\Catsub{#1}{D}{\catba}}
\newcommand{\Dbb}[1][R]{\Catsub{#1}{D}{\catbb}}
\newcommand{\Db}[1][R]{\Catsub{#1}{D}{\catb}}
\newcommand{\Df}[1][R]{\Catsup{#1}{D}{f}}
\newcommand{\Dfba}[1][R]{\Catsupsub{#1}{D}{f}{\catba}}
\newcommand{\Dfbb}[1][R]{\Catsupsub{#1}{D}{f}{\catbb}}
\newcommand{\Dfb}[1][R]{\Catsupsub{#1}{D}{f}{\catb}}
\newcommand{\Dm}[1][R]{\Catsup{#1}{D}{$\fm$}}
\newcommand{\Dmtor}[1][R]{\Catsup{#1}{D}{$\fm$-tor}}
\newcommand{\Dmcpl}[1][R]{\Catsup{#1}{D}{$\fm$-com}}
\newcommand{\Dl}[1][R]{\Catsup{#1}{D}{$\ell$}}
\newcommand{\Dlba}[1][R]{\Catsupsub{#1}{D}{$\ell$}{\catba}}
\newcommand{\Dlbb}[1][R]{\Catsupsub{#1}{D}{$\ell$}{\catbb}}

\newcommand{\Dart}[1][R]{\Catsup{#1}{D}{art}}

\newcommand{\Dartbb}[1][R]{\Catsupsub{#1}{D}{art}{\catbb}}

\newcommand{\km}{\boldsymbol{k}}
\newcommand{\wfm}{\widehat{\fm}}
\newcommand{\wK}{\widehat{K}}
\newcommand{\wN}{\widehat{N}}
\newcommand{\dcmca}[1]{\cite[#1]{dcmca}}
\newcommand{\Dmhat}{\Catsup{\Rhat}{D}{$\wfm$}}
\title{Amplitude inequalities for local (co)homology}

\author[L.W.\ Christensen]{Lars Winther Christensen} %
\address{L.W.C. \ Texas Tech University, Lubbock, TX 79409, U.S.A.}
\email{lars.w.christensen@ttu.edu}
\urladdr{http://www.math.ttu.edu/~lchriste}

\author[L.\ Ferraro]{Luigi Ferraro} %
\address{L.F. \ University of Texas Rio Grande Valley, Edinburg, TX 78539, U.S.A.}
\email{luigi.ferraro@utrgv.edu}
\urladdr{https://faculty.utrgv.edu/luigi.ferraro}

\author[H.\ Holm]{Henrik Holm} %
\address{H.H. \ University of Copenhagen, Universitetsparken 5,
2100 Copenhagen {\O}, Denmark}
\email{holm@math.ku.dk} %
\urladdr{https://web.math.ku.dk/~holm/}

\thanks{L.W.C.\ was partly supported by Simons Foundation
  collaboration grant 428308.}

\thanks{L.F. was partly supported by the Simons Foundation grant MPS-TSM-00007849}

\thanks{H.H. was partly supported by the Danish National Research Foundation through the Copenhagen Centre for Geometry and Topology (DNRF151).}

\date{1 June 2026}



\begin{document}

\begin{abstract}
  Peskine and Szpiro's Intersection Theorem for finitely generated
  modules was generalized by Foxby to all modules and bounded
  complexes of such. We strengthen Foxby's result and prove a dual
  result, which for complexes with finitely generated homology
  recovers Iversen's Amplitude Inequality but also applies to derived
  complete complexes.
\end{abstract}

\maketitle


\section*{Introduction}

\noindent
Let $R$ be a commutative noetherian local ring.  The Intersection
Theorem says that for finitely generated $R$-modules $M$ and $N$ there
is an inequality,
\begin{equation*}
  \pd{M} + \dim{\tpp{N}{M}} \dge \dim{N} \:.
\end{equation*}
It was first proved by Peskine and Szpiro \cite{CPsLSz73} in the
equicharacteristic case and later by Roberts \cite{PRb87,rob} in the
general case. Foxby \cite{HBF79} proved that the inequality,
\begin{equation}
  \label{eq:HBF}
  \fd{M} + \dim{\Ltpp{N}{M}} \dge \dim{N} \:,
\end{equation}
holds for $R$-complexes $M$ and $N$ with bounded homology, where only
the homology of $N$ is assumed to be finitely generated. Applied to
finitely generated $R$-modules, Foxby's version recovers the original
Intersection Theorem.

The Krull dimension of a complex with degreewise finitely generated
homology is captured by vanishing of local cohomology, but in general
only the inequality
\begin{equation*}
  \dim{M} \dge \sup\setof{n \in \ZZ}{\lc[\fm]{n}{M} \ne 0} 
\end{equation*}
holds, and it may be strict. In this paper, we revisit Foxby's work,
switching focus from the Krull dimension to the supremum of
nonvanishing local cohomology. This allows us to improve Foxby's
result, see \thmref{amp-RGam-Ltp} and \rmkref{HBF}; in particular, one
can in \eqref{HBF} replace the Krull dimension by the supremum of
nonvanishing local cohomology.

Even better, the focus on local cohomology suggests dual results
related to nonvanishing of local homology. One such result,
\thmref{amp-LLam-Ltp}, is an inequality that, when specialized to
complexes with degreewise finitely generated homology, recovers
Iversen's Amplitude Inequality~\cite{BIv77} but also applies to
derived complete complexes. To be specific, let $M$ be a complex with
bounded homology and finite flat dimension; \corref{amp-LLam-Ltp-12}
says that if $M$ is derived $\fm$-complete or has degreewise finitely
generated homology, then the inequality
\begin{equation*}
  \amp{\Ltpp{M}{N}} \dge \amp{N}
\end{equation*}
holds for complexes $N$ with degreewise finitely generated homology.

In another direction we prove in \corref{amp-LLam} that derived
complete modules of finite injective dimension only can be found over
Cohen--Macaulay rings. This supplements the fact, first suggested by
Bass, that finitely generated modules of finite injective dimension
only exist over Cohen--Macaulay rings.
\begin{equation*}
  \ast \ \ \ast \ \ \ast
\end{equation*}

\noindent
Throughout this paper, $R$ is a commutative noetherian local ring with
maximal ideal $\fm$ and residue field $\km = R/\fm$. The $\fm$-adic
completion of $R$ is denoted $\Rhat$; it is a local ring with maximal
ideal $\wfm = \fm\Rhat$ and residue field $\Rhat/\wfm \is \km$. We
work in the derived category, $\D$, over $R$ and use homological
notation for $R$-complexes; we refer the reader to \cite{dcmca} for
any unexplained notation and terminology.

It is a well-established fact in local algebra that homological
properties of finitely generated $R$-modules are detected by
(non)vanishing of (co)homology with coefficients in $\km$. For
example, a finitely generated $R$-module $M$ with
$\Tor{1}{\km}{M} = 0$ is free. In this paper we focus on a wide class
of $R$-complexes, one that strictly includes complexes with degreewise
finitely generated homology, and yet enjoys close homological
resemblance to finitely generated modules, see \pgref{Dm} below. The
importance of these complexes was already established in \cite{HBF79}.

The left derived $\fm$-completion functor is denoted $\LLam[\fm]{}$,
and an $R$-complex $M$ with $\LLam[\fm]{M}$ canonically isomorphic to
$M$ is called \emph{derived $\fm$-complete}. Similarly, the right
derived $\fm$-torsion functor is denoted $\RGam[\fm]{}$, and an
$R$-complex $M$ with $\RGam[\fm]{M}$ canonically isomorphic to $M$ is
called \emph{derived $\fm$-torsion}. The study of local cohomology,
$\lc[\fm]{n}{M} = \H[-n]{\RGam[\fm]{M}}$, goes back to Grothendieck
and Hartshorne \cite{lch}. The derived torsion and derived completion
functors and their interactions in commutative algebra, and more
broadly in ring theory, was developed by, primarily, Alonso Tarr\'io,
Jerem\'ias L\'opez, and Lipman \cite{AJL-97}, Greenlees and May
\cite{JPGJPM92}, Porta, Shaul and Yekutieli \cite{PSY-14}, Schenzel
\cite{PSc03}, and Simon \cite{AMS90}. We recall a collection of key
isomorphisms from these works; they can also be found in
\dcmca{13.1.18, 13.3.19, 13.4.12, and 13.4.20}.

\begin{ipg}
  \label{RGLL}
  For $R$-complexes $M$ and $N$ there are isomorphisms,
  \begin{equation}
    \label{eq:RG}
    \Ltp{(\RGam[\fm]{N})}{M} \dqis \RGam[\fm]{\Ltpp{N}{M}} \dqis \Ltp{N}{(\RGam[\fm]{M})}
  \end{equation}
  and
  \begin{equation}
    \label{eq:LL}
    \RHom{\RGam[\fm]{N}}{M} \dqis \LLam[\fm]{\RHom{N}{M}} \dqis \RHom{N}{\LLam[\fm]{M}} \:.
  \end{equation}
  If $M$ is derived $\fm$-torsion, then there are isomorphisms,
  \begin{equation}
    \label{eq:GLb}
    \RHom{M}{\RGam[\fm]{N}} \dqis \RHom{M}{N} \dqis \RHom{M}{\LLam[\fm]{N}}
  \end{equation}
  and
  \begin{equation}
    \label{eq:GLc}
    \Ltp{M}{\RGam[\fm]{N}} \dqis \Ltp{M}{N} \dqis \Ltp{M}{\LLam[\fm]{N}} \:.
  \end{equation}
\end{ipg}

The depth is a classic invariant in local algebra. Through its
cohomological characterization it generalizes naturally to complexes,
and it is already present in \cite{HBF79} and developed further by
Foxby and Iyengar in \cite{HBFSIn03}. The dual notion of width made
its first official appearance in Yassemi's \cite{SYs98b}. For the
equalities recalled below we refer to \dcmca{16.2.3, 16.2.5, 16.2.14,
  16.2.16, and 16.2.34}.

\begin{ipg}
  For an $R$-complex $M$ the depth invariant can be computed in
  several ways:
  \begin{equation}
    \label{eq:dpt1}
    \dpt{M} \deq -\sup{\RHom{\km}{M}} \deq -\sup{\RGam[\fm]{M}} \dge -\sup{M}
  \end{equation}
  and further one has
  \begin{equation}
    \label{eq:dpt2}
    \dpt{\RGam[\fm]{M}} \deq \dpt{M} \deq \dpt{\LLam[\fm]{M}} \:.
  \end{equation}
  Finally, the depth yields a bound on nonvanishing of local homology:
  \begin{equation}
    \label{eq:dpt3}
    \sup{\LLam[\fm]{M}} \dle \dimR -\dpt{M} \:.
  \end{equation}
  
  The width invariant can similarly be computed as
  \begin{equation}
    \label{eq:wdt1}
    \wdt{M} \deq \inf{\Ltpp{\km}{M}} \deq \inf{\LLam[\fm]{M}} \dge \inf{M}
  \end{equation}
  and further one has
  \begin{equation}
    \label{eq:wdt2}
    \wdt{\LLam[\fm]{M}} \deq \wdt{M} \deq \wdt{\RGam[\fm]{M}} \:.
  \end{equation}
  Finally, the width yields a bound on nonvanishing of local
  cohomology:
  \begin{equation}
    \label{eq:wdt3}
    -\inf{\RGam[\fm]{M}} \dle \dimR -\wdt{M} \:.
  \end{equation}
\end{ipg}

We recall that the support and cosupport of an $R$-complex $M$ are the
sets
\begin{align*}
  \supp{M} & \deq \setof{\fp\in\SpecR}{\H{\Ltp{\km(\fp)}{M}} \ne 0}
             \quad\text{and} \\
  \cosupp{M} & \deq \setof{\fp\in\SpecR}{\H{\RHom{\km(\fp)}{M}} \ne 0} \:,
\end{align*}
where $\km(\fp) = (R/\fp)_\fp$ is the residue field of $R_\fp$.  In
the context of commutative algebra, the support appeared in
\cite{HBF79} under the name ``small support'' to emphasize that it is
different from, actually a subset of, the classic notion of support
defined in terms of localization.  The importance of the dual notion
of cosupport was recognized, still in the context of commutative
algebra, by Benson, Iyengar, and Krause \cite{BIK-12}. An elementary
property of both support and cosupport is the ability to recognize
acyclic complexes. Indeed, an $R$-complex $M$ is \emph{acyclic} if it
satisfies the following equivalent conditions:
\begin{center}
  \eqclbl{i}\hspace{.2pc} $\supp{M} \deq \varnothing$ \qquad
  \eqclbl{ii}\hspace{.2pc} $\H{M} = 0$\qquad \eqclbl{iii}\hspace{.2pc}
  $\cosupp{M} \deq \varnothing \:.$
\end{center}
This is established in \cite{BIK-12,HBF79}; for convenience we refer
to \dcmca{15.1.15 and 15.2.8}. The facts recalled below can be found
in \dcmca{13.2.5, 16.1.20, 16.1.34, 16.2.5, and 16.2.27}.

\begin{ipg}
  \label{Dm}
  Complexes that satisfy the equivalent conditions below have
  important homological properties in common with finitely generated
  $R$-modules.
  \begin{center}
    \begin{tabular}[c]{rlrl}
      \eqclbl{i}\hspace{.2pc} & $\fm \in \supp{M}$ 
      & \eqclbl{iii}\hspace{.2pc} & $\wdt{M} < \infty$ \\[2pt]
      \eqclbl{ii}\hspace{.2pc} & $\fm \in \cosupp{M}\qquad$
      & \eqclbl{iv}\hspace{.2pc} & $\dpt{M} < \infty$
    \end{tabular}
  \end{center}
  \vspace{.2pc}

  \noindent We denote by $\Dm$ the full subcategory of $\D$ whose
  objects satisfy the equivalent conditions above. Beware: $\Dm$ is
  not triangulated subcategory of $\D$; it is not even additive, as it
  does not contain the zero complex. It contains the non-acyclic
  complexes in the subcategories $\Dmtor$ of derived $\fm$-torsion
  complexes and $\Dmcpl$ of derived $\fm$-complete complexes. They can
  also be characterized in terms of support and cosupport:
  \begin{equation}
    \label{eq:Dmtor}
    \begin{aligned}
      \Dmtor & \deq \setof{M\in\D}{\supp{M} \subseteq \set{\fm}} \quad\text{and} \\
      \Dmcpl & \deq \setof{M\in\D}{\cosupp{M} \subseteq \set{\fm}} \:.
    \end{aligned}
  \end{equation}
  The category $\Dm$ also contains the nonacyclic complexes from the
  subcategory $\Df$ of complexes with degreewise finitely generated
  homology. One has
  \begin{equation}
    \label{eq:Df}
    \LLam[\fm]{M} \dqis \tp{\Rhat}{M} \qand \wdt{M} \deq \inf{M} 
    \quad\text{ for }M \in \Df \:.
  \end{equation}
\end{ipg}

We recall from \dcmca{15.1.16, 15.1.27, 15.2.9, 15.2.18, 16.2.9, and
  16.2.24} some frequently used properties of complexes in $\Dm$.

\begin{ipg}
  Let $M$ and $N$ be complexes in $\Dm$. The complexes
  \begin{equation}
    \label{eq:Dm}
    \LLam[\fm]{M}\,,\quad \RGam[\fm]{M}\,,\quad \RHom{N}{M}\,,\qand \Ltp{N}{M}
  \end{equation}
  also belong to $\Dm$, and there are equalities:
  \begin{equation}
    \label{eq:dpt-RHom}
    \dpt{\RHom{N}{M}} \deq \wdt{N} + \dpt{M}
  \end{equation}
  and
  \begin{equation}
    \label{eq:wdt-Ltp}
    \wdt{\Ltpp{N}{M}} \deq \wdt{N} + \wdt{M} \:.
  \end{equation}
\end{ipg}

\begin{ipg}
  \label{cmd}
  In \cite{HBF79} Foxby defined the Krull dimension of an $R$-complex
  $M$ and characterized it as follows, see also \dcmca{14.2.1}:
  \begin{equation}
    \label{eq:dim}
    \dim{M} \deq \sup\setof{\dim{\H[n]{M}}-n}{n\in\ZZ} \:.
  \end{equation}
  Foxby also extended Grothendieck Vanishing to complexes, see also
  \dcmca{18.3.22}:
  \begin{equation}
    \label{eq:GV}
    \dim{M} \dge -\inf{\RGam[\fm]{M}} \quad\text{with equality for } M \in \Df \:.
  \end{equation}
  The difference between the Krull dimension and the depth,
  \begin{equation}
    \label{eq:cmd}
    \cmd{M} \deq \dim{M} - \dpt{M} \:,
  \end{equation}
  is called the \emph{Cohen--Macaulay defect} of $M$. For $M\in\Dm$
  the Cohen--Macaulay defect is nonnegative, see \dcmca{18.3.31}, and
  \eqref{dpt1} and \eqref{GV} combine to yield
  \begin{equation}
    \label{eq:cmdRG}
    \cmd{M} \dge \amp{\RGam[\fm]{M}} \quad\text{with equality for }
    M \in \Df \:.
  \end{equation}
\end{ipg}

\begin{ipg}
  We use the notation $\Dl$ for the subcategory of complexes with
  degreewise finite length homology, and we recall from
  \dcmca{16.1.46} that one has
  \begin{equation}
    \label{eq:Dlmtor}
    \Dl \,\subseteq\, \Dmcpl \,\cap\, \Dmtor \:.
  \end{equation}
  For a derived $\fm$-torsion complex $M$ one has $\sup{M} = -\dpt{M}$
  by \eqref{dpt1} and
  \begin{equation}
    \label{eq:RGdim}
    \inf{M} \deq - \dim{M}
  \end{equation}
  by \eqref{dim}, see also \dcmca{16.1.34(b)}; thus the next equality
  holds,
  \begin{equation}
    \label{eq:Dl}
    \cmd{M} \deq \amp{M} \:.
  \end{equation}
  By \eqref{Dlmtor} this equality holds, in particular, for $M$ in
  $\Dl$.
\end{ipg}


\section{Depth and width}

\noindent
We start by collecting a smattering of homological properties of
derived $\fm$-torsion complexes and derived $\fm$-complete complexes;
they will be used throughout.

\begin{prp}
  \label{prp:RG}
  Let $M$ be an $R$-complex; the following assertions hold.
  \begin{prt}
  \item The complex $\RGam[\fm]{M}$ is an $\Rhat$-complex and derived
    $\wfm$-torsion.
  \item If $M$ belongs to $\Dm$, then $\RGam[\fm]{M}$ belongs to
    $\Dmhat$.
  \item There are isomorphisms,
    \begin{equation*}
      \Ltp[\Rhat]{\km}{\RGam[\fm]{M}} \dqis \Ltp{\km}{\RGam[\fm]{M}}
      \dqis \Ltp{\km}{M}\\[-.5\baselineskip]
    \end{equation*}
    and
    \begin{equation*}
      \RHom[\Rhat]{\km}{\RGam[\fm]{M}} \dqis \RHom{\km}{\RGam[\fm]{M}}
      \dqis \RHom{\km}{M} \:. 
    \end{equation*}
  \item There are (in)equalities,
    \begin{equation*}
      \sup{\RGam[\fm]{M}} \dle \fd[\Rhat]{\RGam[\fm]{M}} \deq
      \fd{\RGam[\fm]{M}} \deq \sup{\Ltpp{\km}{M}} \:.
    \end{equation*}
  \item If $\Ltp{\km}{M}$ is in $\Db$, then so is $\RGam[\fm]{M}$ and
    $\fd{\RGam[\fm]{M}}$ is finite.
  \item If the complex $\RHom{\km}{M}$ belongs to $\Dba$, then there
    are (in)equalities,
    \begin{equation*}
      -\inf{\RGam[\fm]{M}} \dle \id[\Rhat]{\RGam[\fm]{M}} \deq
      \id{\RGam[\fm]{M}} \deq -\inf{\RHom{\km}{M}} \:.
    \end{equation*}
  \item If\, $\RHom{\km}{M}$ is in $\Db$, then so is $\RGam[\fm]{M}$
    and $\id{\RGam[\fm]{M}}$ is finite.
  \end{prt}
\end{prp}

\begin{prf*}
  The assertions in part (a) hold by \dcmca{11.3.18 and 13.3.23(a)}.

  \proofoftag{c} In the first display, the first isomorphism holds by
  \dcmca{13.4.18(c)}, and as $\km$ per \eqref{Dlmtor} is derived
  $\fm$-torsion, the second isomorphism holds by
  \eqref{GLc}. Similarly, the isomorphisms in second display hold by
  \dcmca{13.4.18(b)} and \eqref{GLb}.

  \proofoftag{b} For $M$ in $\Dm$ the complexes in (c) are not
  acyclic, so $\RGam[\fm]{M}$ is in $\Dmhat$.

  \proofoftag{d} The inequality is standard and the equalities hold by
  \dcmca{15.4.25 and 17.3.7}.

  \proofoftag{e} The assertions follow from part (d) and the
  (in)equalities,
  \begin{equation*}
    \inf{\RGam[\fm]{M}} \dge \wdt{M} - \dimR \deq \inf{\Ltpp{\km}{M}} - \dimR \:,
  \end{equation*}
  which hold by \eqref{wdt3} and \eqref{wdt1}.

  \proofoftag{f} The inequality is standard and the equalities hold in
  view of \eqref{dpt1} by \dcmca{16.1.23 and 16.3.13}.

  \proofoftag{g} The assertions follow from part (f) and \eqref{dpt1}.
\end{prf*}

\begin{prp}
  \label{prp:LL}
  Let $M$ be an $R$-complex; the following assertions hold.
  \begin{prt}
  \item The complex $\LLam[\fm]{M}$ is an $\Rhat$-complex and derived
    $\wfm$-complete.

  \item If $M$ belongs to $\Dm$, then $\LLam[\fm]{M}$ belongs to
    $\Dmhat$.
    
  \item There are isomorphisms,
    \begin{equation*}
      \Ltp[\Rhat]{\km}{\LLam[\fm]{M}} \dqis \Ltp{\km}{\LLam[\fm]{M}}
      \dqis \Ltp{\km}{M}\\[-.5\baselineskip]
    \end{equation*}
    and
    \begin{equation*}
      \RHom[\Rhat]{\km}{\LLam[\fm]{M}} \dqis \RHom{\km}{\LLam[\fm]{M}}
      \dqis \RHom{\km}{M} \:.
    \end{equation*}

  \item There are (in)equalities,
    \begin{equation*}
      \hspace{5ex}
      -\inf{\LLam[\fm]{M}} \dle \id[\Rhat]{\LLam[\fm]{M}} \deq
      \id{\LLam[\fm]{M}} \deq -\inf{\RHom{\km}{M}} \:.
    \end{equation*}

  \item If\, $\RHom{\km}{M}$ is in $\Db$, then so is $\LLam[\fm]{M}$
    and $\id{\LLam[\fm]{M}}$ is finite.

  \item If $\Ltp{\km}{M}$ belongs to $\Dbb$, then there are
    (in)equalities,
    \begin{equation*}
      \hspace{5ex}
      \sup{\LLam[\fm]{M}} \dle \fd[\Rhat]{\LLam[\fm]{M}} \deq \fd{\LLam[\fm]{M}} \deq
      \sup{\Ltpp{\km}{M}} \:.
    \end{equation*}

  \item If $\Ltp{\km}{M}$ is in $\Db$, then so is $\LLam[\fm]{M}$ and
    $\fd{\LLam[\fm]{M}}$ is finite.
  \end{prt}
\end{prp}

\begin{prf*}
  The assertions in part (a) hold by \dcmca{11.3.4 and 13.1.21(a)}.

  \proofoftag{c} In the first display, the first isomorphism holds by
  \dcmca{13.4.18(c)}, and as $\km$ per \eqref{Dlmtor} is derived
  $\fm$-torsion, the second isomorphism holds by
  \eqref{GLc}. Similarly, the isomorphisms in second display hold by
  \dcmca{13.4.18(b)} and \eqref{GLb}.
  
  \proofoftag{b} For $M$ in $\Dm$ the complexes in (c) are not
  acyclic, so $\LLam[\fm]{M}$ is in $\Dmhat$.

  \proofoftag{d} The inequality is standard and the equalities hold by
  \dcmca{17.3.20 and 17.3.23}.

  \proofoftag{e} The assertions follow from part (d) and the
  (in)equalities,
  \begin{equation*}
    \sup{\LLam[\fm]{M}} \dle \dimR - \dpt{M} \deq \dimR + \sup{\RHom{\km}{M}} \:,
  \end{equation*}
  which hold by \eqref{dpt3} and \eqref{dpt1}.

  \proofoftag{f} The inequality is standard and the equalities hold in
  view of \eqref{wdt1} by \dcmca{16.1.24 and 16.3.18}.

  \proofoftag{g} The assertions follow from part (f) and \eqref{wdt1}.
\end{prf*}

In the balance of this section, we establish some formulas for
computing depth and width of derived Hom and tensor product complexes.
The first result below extends \dcmca{16.3.1(a)} and its corollary
\dcmca{16.3.3}. In the terminology of Frankild and J{\o}rgensen
\cite{AFrPJr03b}, the assumption on $M$ is that it is a complex of
finite $\km$-projective dimension. Indeed, by Matlis duality one has
\begin{equation*}
  -\inf{\RHom{M}{\km}} \deq \sup{\Ltpp{\km}{M}}\:,
\end{equation*} see \dcmca{16.3.8}, and that quantity is assumed to be finite.

\begin{prp}
  \label{prp:ABE}
  Let $M$ be an $R$-complex with $\Ltp{\km}{M}$ in $\Dba$. One has
  \begin{equation}
    \tag{a}
    \dpt{M} \deq \dptR - \sup{\Ltpp{\km}{M}} \,>\, -\infty\:.
  \end{equation}
  Moreover, for every $R$-complex $N$ with $\RHom{\km}{N}$ in $\Dba$
  there is an equality,
  \begin{equation}
    \tag{b}
    \dpt{\Ltpp{N}{M}} \deq \dpt{N} + \dpt{M} - \dptR \:.
  \end{equation}
\end{prp}

\begin{prf*}
  The $R$-module $\km$ is derived $\fm$-torsion, see \eqref{Dlmtor},
  so the first equality below follows from \eqref{GLc}. By
  \prpref{RG}(d) the $R$-complex $\RGam[\fm]{M}$ has finite flat
  dimension, whence \dcmca{16.3.3} yields the second equality. The
  third and final equality holds by \eqref{dpt2}.
  \begin{align*}
    \sup{\Ltpp{\km}{M}} 
    & \deq \sup{\Ltpp{\km}{\RGam[\fm]{M}}} \\
    & \deq \dptR - \dpt{\RGam[\fm]{M}} \\
    & \deq \dptR - \dpt{M} \:.
  \end{align*}
  Rearranging the terms yields the equality asserted in part (a), and
  the inequality holds by the assumption on $M$.
  
  For part (b), notice that $\dpt{N} > -\infty$ holds by
  \eqref{dpt1}. Now the asserted equality follows from the next
  computation where the equalities hold, in order, by \eqref{dpt2},
  \eqref{RG}, \dcmca{16.3.1(a)}, and \eqref{dpt2}.
  \begin{align*}
    \dpt\Ltpp{N}{M} 
    & \deq\dpt\RGam[\fm]{\Ltpp{N}{M}} \\
    & \deq\dpt\Ltpp{N}{\RGam[\fm]M} \\
    & \deq\dpt N+\dpt\RGam[\fm]M-\dptR \\
    & \deq\dpt N+\dpt M-\dptR\:. \qedhere
  \end{align*}
\end{prf*}

Every complex $M$ in $\Dm$ satisfies the inequality
\begin{equation*}
  \dpt{M} + \wdt{M} \dle \dimR \:,
\end{equation*}
which Strooker \corcite[6.1.10]{JRS90} credits to Bartijn; see also
Schenzel and Simon \prpcite[5.3.9]{PScAMS} or \dcmca{16.2.35}. The
next result and its partner, \corref{coABE}, improve the bound under
additional assumptions on $M$.

\begin{cor}
  \label{cor:ABE}
  Let $M$ be an $R$-complex. If $\Ltp{\km}{M}$ is in $\Dba$, then one
  has
  \begin{equation*}
    \dpt{M} + \wdt{M} \deq \dptR - \amp{\Ltpp{\km}{M}} \:.
  \end{equation*}
\end{cor}

\begin{prf*}
  One has $\wdt{M} = \inf{\Ltpp{\km}{M}}$, see \eqref{wdt1}, and
  adding this quantity to the equality in \prpref{ABE}(a) yields the
  asserted equality.
\end{prf*}

\begin{rmk}
  \label{rmk:HBF}
  In \cite{hha} Foxby showed that in the inequality \eqref{HBF} one
  can replace $\fd{M}$ with the smaller quantity $\dptR - \dpt{M}$;
  see \dcmca{16.3.4}. Doing so and subtracting the equality from
  \prpref{ABE}, on the form
  \begin{equation*}
    \dptR - \dpt{M} + \dpt{\Ltpp{N}{M}} \deq \dpt{N} \:,
  \end{equation*}
  one arrives at
  \begin{equation*}
    \cmd{\Ltpp{N}{M}}  \dge \cmd{N} \:,
  \end{equation*} 
  which is how Foxby's version of the Intersection Theorem is stated
  in \cite{dcmca,hha}. It serves as a template for the main results in
  this paper, including \thmref{amp-RGam-Ltp}, which improves Foxby's
  version by replacing $\cmd{\Ltpp{N}{M}}$ by
  $\amp{\RGam[\fm]{\Ltpp{N}{M}}}$, cf.~\eqref{cmd}, and by adding an
  upper bound on this quantity based on \corref{ABE}.
\end{rmk}

The next result extends \dcmca{16.3.9(a)} and its corollary
\dcmca{16.3.11}. The condition on $M$ is, in the terminology of
\cite{AFrPJr03b}, that it is a complex of finite $\km$-injective
dimension.
\begin{prp}
  \label{prp:coABE}
  Let $M$ be an $R$-complex with $\RHom{\km}{M}$ in $\Dbb$. One~has
  \begin{equation}
    \tag{a}
    \wdt M \deq \dptR + \inf{\RHom{\km}{M}} \,>\, -\infty\:.
  \end{equation}
  Moreover, for every $R$-complex $K$ with $\RHom{\km}{K}$ in $\Dba$
  there is an equality,
  \begin{equation}
    \tag{b}
    \wdt{\RHom{K}{M}} \deq \dpt{K} + \wdt{M} - \dptR \:.
  \end{equation}
\end{prp}

\begin{prf*}
  The $R$-module $\km$ is derived $\fm$-torsion, see \eqref{Dlmtor},
  so the first equality below follows from \eqref{GLb}. The
  $R$-complex $\LLam[\fm]{M}$ has finite injective dimension by
  \prpref{LL}(d), and hence \dcmca{16.3.11} yields the second
  equality. The third and final equality holds by \eqref{wdt2}.
  \begin{align*}
    -\inf{\RHom{\km}{M}} 
    & \deq -\inf{\RHom{\km}{\LLam[\fm]{M}}} \\
    & \deq \dptR - \wdt{\LLam[\fm]{M}} \\
    & \deq \dptR - \wdt{M} \:.
  \end{align*}
  Rearranging the terms yields the equality asserted in part (a), and
  the inequality holds by the assumption on $M$.

  For part (b), notice that $\dpt{K}>-\infty$ holds by
  \eqref{dpt1}. Now the asserted equality follows from the next
  computation where the equalities hold, in order, by \eqref{wdt2},
  \eqref{LL}, \dcmca{16.3.9(a)}, and \eqref{wdt2}.
  \begin{align*}
    \wdt{\RHom{K}{M}} 
    & \deq \wdt{\LLam[\fm]{\RHom{K}{M}}} \\
    & \deq \wdt{\RHom{K}{\LLam[\fm]{M}}} \\
    & \deq \dpt{K} + \wdt{\LLam[\fm]{M}} - \dptR \\
    & \deq \dpt{K} + \wdt{M} - \dptR \:. \qedhere
  \end{align*}
\end{prf*}

\begin{cor}
  \label{cor:coABE}
  Let $M$ be an $R$-complex. If\, $\RHom{\km}{M}$ is in $\Dbb$, then
  one has
  \begin{equation*}
    \dpt{M} + \wdt{M} \deq \dptR - \amp{\RHom{\km}{M}} \:.
  \end{equation*}
\end{cor}

\begin{prf*}
  One has $\dpt{M} = -\sup{\RHom{\km}{M}}$, see \eqref{dpt1}, and
  adding this quantity to the equality in \prpref{coABE}(a) yields the
  asserted equality.
\end{prf*}

The next result, which extends \dcmca{16.3.5(a)}, supplements parts
(b) in \prpref[Propositions~]{ABE} and \prpref[]{coABE}. There is no
part (a) stated, as it would only duplicate part (a) in \prpref{ABE},
and for
the same reason no corollary is stated; it would just rephrase
\corref{ABE}.

\begin{prp}
  \label{prp:coABE-bis}
  Let $M$ be an $R$-complex such that $\Ltp{\km}{M}$ belongs to
  $\Dba$. For every $R$-complex $K$ with $\Ltp{\km}{K}$ in $\Dbb$
  there is an equality,
  \begin{equation*}
    \wdt{\RHom{M}{K}} \deq \dpt{M} + \wdt{K} - \dptR \:.
  \end{equation*}
\end{prp}

\begin{prf*}
  By \prpref{RG}(d) the complex $\RGam[\fm]{M}$ has finite flat
  dimension and, thus, finite projective dimension as $R$ is
  local. Further, $\wdt{K}>-\infty$ holds by \eqref{wdt1}. Now the
  asserted equality follows from the next computation where the
  equalities hold, in order, by \eqref{wdt2}, \eqref{LL},
  \dcmca{16.3.5(a)}, and \eqref{dpt2}.
  \begin{align*}
    \wdt{\RHom{M}{K}} 
    & \deq \wdt{\LLam[\fm]{\RHom{M}{K}}} \\
    & \deq \wdt{\RHom{\RGam[\fm]{M}}{K}} \\
    & \deq \dpt{\RGam[\fm]{M}} + \wdt{K} - \dptR \\
    & \deq \dpt{M} + \wdt{K}  - \dptR \:. \qedhere
  \end{align*}
\end{prf*}


\section{Nonvanishing of local homology}

\noindent
For complexes in $\Df$ the next result is, in view of \eqref{Df},
known from \dcmca{18.5.7}.

\begin{prp}
  \label{prp:amp-LLam}
  Let $M$ a complex in $\Dm$. If the complex $\RHom{\km}{M}$ belongs
  to $\Dbb$, then there is an inequality,
  \begin{equation*}
    \amp{\LLam[\fm]{M}} \dge \cmdR \:.
  \end{equation*}
\end{prp}

\begin{prf*}
  Let $H$ be an $R$-module of maximal depth, i.e.\ $\dpt{H} = \dimR$;
  such modules exist by work of Andr\'e \cite{YAn18}, see
  \dcmca{18.4.9}. The complex $\RHom{H}{\LLam[\fm]{M}}$ belongs to
  $\Dm$, see \eqref{Dm}; in particular, it is not acyclic which
  explains the second inequality in the computation below. The first
  inequality holds by \dcmca{7.6.7}. The first equality is trivial,
  and the second follows from \eqref{LL}.  The fourth equality comes
  from \prpref{coABE}, and the penultimate equality holds by
  \eqref{cmd}. The remaining equalities hold by \eqref{wdt1}.
  \begin{align*}
    \sup{\LLam[\fm]{M}} 
    & \deq \sup{\LLam[\fm]{M}} - \inf{H} \\
    & \dge \sup{\RHom{H}{\LLam[\fm]{M}}} \\
    & \dge \inf{\RHom{H}{\LLam[\fm]{M}}} \\
    & \deq \inf{\LLam[\fm]{\RHom{H}{M}}} \\
    & \deq \wdt{\RHom{H}{M}} \\
    & \deq \dpt{H} + \wdt{M} - \dptR \\
    & \deq \cmdR + \wdt{M} \\
    & \deq \cmdR + \inf{\LLam[\fm]{M}} \:.
  \end{align*}
  The asserted inequality now follows by rearranging the terms.
\end{prf*}

It was suggested by Bass \cite{HBs63} and proved by Peskine and Szpiro
\cite{CPsLSz73} in the equicharacteristic case, and in general by
Roberts \cite{PRb87, rob}, that existence of a nonzero finitely
generated $R$-module of finite injective dimension forces $R$ to be
Cohen--Macaulay. The next corollary elaborates on this classic
result. For perspective, we recall from \dcmca{13.1.33} that an
$\fm$-complete module is derived $\fm$-complete.

\begin{cor}
  \label{cor:amp-LLam}
  If there exists a nonzero derived $\fm$-complete or finitely
  generated $R$-module of finite injective dimension, then $R$ is
  Cohen--Macaulay.
\end{cor}

\begin{prf*}
  Let $M \ne 0$ be an $R$-module of finite injective dimension. If $M$
  is derived $\fm$-complete, then one has
  $\amp{\LLam[\fm]{M}} = \amp{M} = 0$.  If $M$ is finitely generated,
  then there is an isomorphism $\LLam[\fm]{M} \qis \tp{\Rhat}{M}$, see
  \eqref{Df}, so again $\amp{\LLam[\fm]{M}} = 0$ holds. Now invoke
  \prpref{amp-LLam}.
\end{prf*}

The proof of the next theorem is inspired by Foxby's proof of
\thmcite[4.1]{HBF79}.

\begin{thm}
  \label{thm:amp-LLam-RHom}
  Let $M$ and $K$ be complexes with $M \in \Dm$ and $K \in \Dfba$. If
  the complex $\RHom{\km}{M}$ belongs to $\Dbb$, then there are
  inequalities,
  \begin{equation*}
    \amp{\RHom{\km}{M}} + \cmd{K}
    \dge \amp{\LLam[\fm]{\RHom{K}{M}}} \dge \cmd{K} \:.
  \end{equation*}
\end{thm}

\begin{prf*}
  Both inequalities are trivial if $K$ is acyclic, so assume that is
  not the case. That is, assume that $K$ belongs to $\Dm$, see
  \ref{Dm}.  As $K$ belongs to $\Dba$ so does $\RHom{\km}{K}$, see
  \dcmca{7.6.7}. Now \eqref{wdt1} and \prpref{coABE}~yield
  \begin{align*}
    -\inf\LLam[\fm]{\RHom{K}{M}} 
    & \deq -\wdt\RHom{K}{M} \\
    & \deq  -\dpt{K} - \wdt{M} + \dptR \:.
  \end{align*}
  Thus, it follows from \corref{coABE} and \eqref{cmd} that it is
  sufficient to establish the inequalities
  \begin{equation}
    \tag{$\ast$}
    \begin{aligned}
      \dim{K} - \dpt{M} 
      & \dge \sup{\LLam[\fm]{\RHom{K}{M}}} \\
      & \dge \dim{K} + \wdt{M} - \dptR \:.
    \end{aligned}
  \end{equation}
  By \eqref{LL} and \eqref{GLb} there are isomorphisms
  \begin{equation*}
    \LLam[\fm]{\RHom{K}{M}} \dqis \RHom{\RGam[\fm]{K}}{M}
    \dqis \RHom{\RGam[\fm]{K}}{\RGam[\fm]{M}} \:.
  \end{equation*}
  The left-hand inequality in $(\ast)$ now follows from \eqref{GV},
  \eqref{dpt1}, and \dcmca{7.6.7}:
  \begin{equation*}
    \dim{K} - \dpt{M} \deq  - \inf{\RGam[\fm]{K}} + \sup{\RGam[\fm]{M}}
    \dge \sup{\LLam[\fm]{\RHom{K}{M}}} \:.
  \end{equation*}
  In the next computation, the first inequality holds by
  \eqref{dpt1}. The equalities follow from \eqref{dpt2} and
  \eqref{dpt-RHom}. The final inequality holds by \eqref{dpt3}.
  \begin{align*}
    \sup{\LLam[\fm]{\RHom{K}{M}}} 
    & \dge -\dpt{\LLam[\fm]{\RHom{K}{M}}} \\
    & \deq -\dpt{\RHom{K}{M}} \\    
    & \deq - \wdt{K} - \dpt{M} \\    
    & \dge - \wdt{K} + \sup{\LLam[\fm]{M}} - \dimR \:.
  \end{align*}
  Per \eqref{Df} one has $\wdt{K} = \inf{K}$, and the computation
  above shows that the right-hand inequality in $(\ast)$ holds if one
  has $\inf{K} = -\infty$ or $\sup{\LLam[\fm]{M}} = \infty$. One can
  henceforth assume that $K$ belongs to $\Dfb$ and $\LLam[\fm]{M}$ to
  $\Db$ as $\inf{\LLam[\fm]{M}} = \wdt{M} > -\infty$ holds by
  \eqref{wdt1} and \prpref{coABE}.
  
  To prove the right-hand inequality in $(\ast)$ we first assume that
  $K$ is a cyclic module $R/\fa$. The $R/\fa$-complex
  $X = \RHom{K}{M}$ meets the assumptions in \prpref{amp-LLam}; indeed
  $\km$ is also the residue field of the local ring $R/\fa$ and one
  has
  \begin{equation*}
    \RHom[R/\fa]{\km}{X} \dqis \RHom{\km}{M} \:.
  \end{equation*}
  As $M$ is in $\Dm$, the complex $\RHom{\km}{M}$ is not acyclic, so
  $X$ is in $\Catsup{R/\fa}{D}{$\fm/\fa$}$ and $\RHom[R/\fa]{\km}{X}$
  belongs to $\Dbb[R/\fa]$. By \eqref{dpt1} one has
  \begin{equation*}
    \sup{\RHom{\km}{K}} \deq  - \dpt{R/\fa} \,<\, \infty \:,
  \end{equation*}
  so the $R$-complexes $K$ and $M$ satisfy the assumptions in
  \prpref{coABE}. In the computation below the first equality is
  trivial and the second holds by \eqref{wdt1} and independence of
  base for local homology, see \dcmca{13.1.21(a)}. The inequality
  follows from \prpref{amp-LLam} and the penultimate equality from
  \prpref{coABE}. The final equality follows from \eqref{cmd}.
  \begin{align*}
    \sup{\LLam[\fm]{\RHom{K}{M}}}
    & \deq \amp{\LLam[\fm]{\RHom{K}{M}}} + \inf{\LLam[\fm]{\RHom{K}{M}}} \\
    & \deq \amp{\LLam[\fm/\fa]{X}} + \wdt{\RHom{K}{M}} \\    
    & \dge \cmdR/\fa + \wdt{\RHom{K}{M}} \\
    & \deq \cmd{K} + \dpt{K} + \wdt{M} - \dptR\\
    & \deq \dim{K}  + \wdt{M} - \dptR \:.
  \end{align*}
  Next we prove the right-hand inequality in $(\ast)$ under the
  assumption that $K$ is a finitely generated $R$-module. Set
  $\fa = (0:_RK)$; the second equality below holds by
  \dcmca{17.6.13}. The first and third equalities follow from
  \eqref{LL}. The inequality was established above and the last
  equality is standard, see \dcmca{14.1.1}.
  \begin{align*}
    \sup{\LLam[\fm]{\RHom{K}{M}}} 
    & \deq \sup{\RHom{K}{\LLam[\fm]{M}}} \\
    & \deq \sup{\RHom{R/\fa}{\LLam[\fm]{M}}} \\
    & \deq \sup{\LLam[\fm]{\RHom{R/\fa}{M}}} \\
    & \dge \dim{R/\fa}  + \wdt{M} - \dptR \\
    & \deq \dim{K} + \wdt{M} - \dptR \:.
  \end{align*}
  Finally we prove the right-hand inequality in $(\ast)$ for any
  complex $K$ in $\Dfb$. The first and third equalities below follow
  from \eqref{LL}. The second equality holds by \dcmca{17.6.14}. The
  inequality was established above, and the final equality is
  \eqref{dim}.
  \begin{multline*}
    \sup{\LLam[\fm]{\RHom{K}{M}}} \\
    \begin{aligned}
      & \deq \sup{\RHom{K}{\LLam[\fm]{M}}} \\
      & \deq \sup\setof{\sup{\RHom{\H[n]{K}}{\LLam[\fm]{M}}} - n}{n\in\ZZ} \\
      & \deq \sup\setof{\sup{\LLam[\fm]{\RHom{\H[n]{K}}{M}}} - n}{n\in\ZZ} \\
      & \dge \sup\setof{\dim{\H[n]{K} - n}}{n\in\ZZ} + \wdt{M} - \dptR \\
      & \deq \dim{K} + \wdt{M} - \dptR \:.
    \end{aligned}\\[-1.05\baselineskip]
    \qedhere
  \end{multline*}
\end{prf*}

The next result has a counterpart in \corref{amp-RGam-RHom-12}.

\begin{cor}
  \label{cor:amp-LLam-RHom-12}
  Let $M$ and $K$ be $R$-complexes with $K \in \Dfba$. Assume that $M$
  is of finite injective dimension and not acyclic. If $M$ is derived
  $\fm$-complete or belongs to $\Dfb$, then there are inequalities,
  \begin{equation*}
    \id{M} - \dpt{M} + \cmd{K} \dge \amp{\RHom{K}{M}} \dge \cmd{K} \:.
  \end{equation*}
\end{cor}

\begin{prf*}
  The complex $M$ belongs to $\Dm$, see \ref{Dm}. If $M$ is derived
  $\fm$-complete, then \eqref{dpt1} and \prpref{LL}(d) yield
  \begin{equation}
    \tag{$\ast$}
    \sup{\RHom{\km}{M}} \deq - \dpt{M} \qand
    \inf{\RHom{\km}{M}} \deq - \id{M} \:;
  \end{equation}
  in particular, the complex $\RHom{\km}{M}$ belongs to
  $\Dbb$. Further, the complex $\RHom{K}{M}$ is derived
  $\fm$-complete, see \eqref{LL}, so the asserted inequalities follow
  from \thmref{amp-LLam-RHom} and the equalities $(\ast)$.

  If $M$ belongs to $\Dfb$, then the equalities $(\ast)$ hold by
  \eqref{dpt1} and \dcmca{16.4.8}. The complex $\RHom{K}{M}$ is in
  $\Dfbb$, see \dcmca{15.4.9}, so there are by \eqref{Df}, and
  faithful flatness of $\Rhat$ as an $R$-module, equalities
  \begin{equation*}
    \amp{\LLam[\fm]{\RHom{K}{M}}} \deq \amp{\tpp{\Rhat}{\RHom{K}{M}}}
    \deq \amp{\RHom{K}{M}} \:.
  \end{equation*}
  Now the asserted inequalities follow from \thmref{amp-LLam-RHom} and
  the equalities $(\ast)$.
\end{prf*}

\begin{exa}
  \label{exa:Dhat}
  The complex $D = \LLam[\fm]{\E{\km}}$ is dualizing for $\Rhat$, see
  \dcmca{18.2.9}. The equalities
  $\id{D} = \id[\Rhat]{D} = \dpt[\Rhat]{D} = \dpt{D}$ hold by
  \dcmca{17.3.20, 18.2.11, and 18.3.1}, so for a complex $K$ in
  $\Dfba$ \corref{amp-LLam-RHom-12} yields
  \begin{equation*}
    \amp{\RHom{K}{D}} \deq \cmd{K} \:.
  \end{equation*}
  Further, by \eqref{LL}, \dcmca{16.1.25(6)}, and \eqref{dpt1} one has
  \begin{equation*}
    \inf{\RHom{K}{D}} \deq \inf{\RHom{\RGam[\fm]{K}}{\E{\km}}} = -\sup{\RGam[\fm]{K}} \deq \dpt{K}
  \end{equation*}
  and, therefore, per \eqref{cmd} also
  \begin{equation*}
    \sup{\RHom{K}{D}} \deq \dim{K} \:.
  \end{equation*}
  These interpretations of supremum and infimum of $\RHom{K}{D}$ are
  standard if one replaces $D$ by a dualizing $R$-complex, see
  \dcmca{18.2.31}, and like in the proof of that statement one can
  verify that they hold for all complexes $K$ in $\Df$.
\end{exa}

The inequalities in the next result also hold under slightly different
conditions on the complexes, see \corref{amp-RGam-RHom-3}.

\begin{cor}
  \label{cor:amp-LLam-RHom-3}
  Let $M$ and $K$ be $R$-complexes with $M \in \Dm$ and $K \in \Dlba$.
  If the complex $\RHom{\km}{M}$ belongs to $\Dbb$, then there are
  inequalities,
  \begin{equation*}
    \amp{\RHom{\km}{M}} + \amp{K} \dge \amp{\RHom{K}{M}} \dge \amp{K} \:.
  \end{equation*}
\end{cor}

\begin{prf*}
  Recall from \eqref{Dlmtor} that $K$ is derived $\fm$-torsion, so
  $\RHom{K}{M}$ is per \eqref{LL} derived $\fm$-complete. Now apply
  \thmref{amp-LLam-RHom} and~\eqref{Dl}.
\end{prf*}

The next result extends \dcmca{13.1.19(c)} for the ideal
$\mathfrak{a} = \fm$.

\begin{prp}
  \label{prp:LLam}
  Let $M$ and $N$ be $R$-complexes with $N \in \Df$. If the complex
  $\Ltp{\km}{M}$ belongs to $\Db$, then there is an isomorphism,
  \begin{equation*}
    \LLam[\fm]{\Ltpp{N}{M}} \dqis \Ltp{N}{\LLam[\fm]{M}} \:.
  \end{equation*}
\end{prp}

\begin{prf*}
  The first two isomorphisms in the display below follow from
  \dcmca{13.4.1(c)} and \eqref{RG}. By \prpref{RG}(e) the $R$-complex
  $\RGam[\fm]{M}$ has bounded homology and finite flat dimension.  As
  $N$ is in $\Df$ the third isomorphism now holds by \dcmca{13.1.19(c)
    and 13.4.1(c)}.
  \begin{align*}
    \LLam[\fm]{\Ltpp{N}{M}} 
    & \dqis \LLam[\fm]{\RGam[\fm]{\Ltpp{N}{M}}} \\
    & \dqis \LLam[\fm]{\Ltpp{N}{\RGam[\fm]{M}}} \\
    & \dqis \Ltp{N}{\LLam[\fm]{M}} \:. \qedhere
  \end{align*}
\end{prf*}

We derive the next result from \thmref{amp-LLam-RHom} by way of
Grothendieck Duality.

\begin{thm}
  \label{thm:amp-LLam-Ltp}
  Let $M$ and $N$ be complexes with $M \in \Dm$ and $N \in \Df$. If
  the complex $\Ltp{\km}{M}$ belongs to $\Db$, then there are
  inequalities,
  \begin{equation*}
    \amp{\Ltpp{\km}{M}} + \amp{N} 
    \dge \amp{\LLam[\fm]{\Ltpp{N}{M}}} \\
    \dge \amp{N} \:.
  \end{equation*}
\end{thm}

\begin{prf*}
  Both inequalities are trivial if $N$ is acyclic, so assume that is
  not the case. That is, assume that $N$ belongs to $\Dm$, see
  \ref{Dm}. The equalities below hold by \eqref{wdt1} and
  \eqref{wdt-Ltp}.
  \begin{equation*}
    \inf{\LLam[\fm]{\Ltpp{N}{M}}}  
    \deq \wdt{\Ltpp{N}{M}} 
    \deq \wdt{N} + \wdt{M} \:.                                               
  \end{equation*}
  Per \eqref{Df} one has $\wdt{N} = \inf{N}$, so the equalities show
  that both of the asserted inequalities are trivial if
  $\inf{N} = -\infty$ holds.

  One can now assume that $N$ belongs to $\Dfbb$.  Set
  $\wN = \tp{\Rhat}{N}$; it is by \dcmca{18.3.2(c)} a complex in
  $\Dfbb[\Rhat]$. In the display below, the first isomorphism holds by
  \prpref{LL}(a), and the second follows from \prpref{LLam}. As
  $\LLam[\fm]{M}$ by another application of \prpref{LL}(a) is a
  complex in $\Dmhat$, the third isomorphism is standard, see
  \dcmca{12.3.28}. The ring $\Rhat$ has a dualizing complex, $D$, and
  the fourth isomorphism is Grothendieck Duality, see
  \dcmca{18.2.3}. The last isomorphism holds by \dcmca{12.3.23(c)}, as
  the $\Rhat$-complex $\LLam[\fm]{M}$ has bounded homology and finite
  flat dimension by \prpref{LL}(f,g).
  \begin{align*}
    \LLam[\fm]{\Ltpp{N}{M}} 
    & \dqis \LLam[\wfm]{\LLam[\fm]{\Ltpp{N}{M}}} \\
    & \dqis \LLam[\wfm]{\Ltpp{N}{\LLam[\fm]{M}}} \\
    & \dqis \LLam[\wfm]{\Ltpp[\Rhat]{\wN}{\LLam[\fm]{M}}} \\
    & \dqis
      \LLam[\wfm]{\Ltpp[\Rhat]{\RHom[\Rhat]{\RHom[\Rhat]{\wN}{D}}{D}}{\LLam[\fm]{M}}} \\
    & \dqis
      \LLam[\wfm]{\RHom[\Rhat]{\RHom[\Rhat]{\wN}{D}}{\Ltp[\Rhat]{D}{\LLam[\fm]{M}}}}
      \:.
  \end{align*}
  The complex $\Ltp[\Rhat]{D}{\LLam[\fm]{M}}$ belongs to $\Dmhat$, see
  \eqref{Dm}.  The first isomorphism below holds by
  \dcmca{12.3.23(c)}. The second isomorphism holds as one can assume
  that $D$ is normalized, see \dcmca{18.2.33 and 18.2.24}, and the
  third comes from \prpref{LL}(c).
  \begin{equation}
    \tag{$\ast$}
    \begin{aligned}
      \RHom[\Rhat]{\km}{\Ltp[\Rhat]{D}{\LLam[\fm]{M}}} 
      & \dqis \Ltp[\Rhat]{\RHom[\Rhat]{\km}{D}}{\LLam[\fm]{M}} \\
      & \dqis \Ltp[\Rhat]{\km}{\LLam[\fm]{M}} \\
      & \dqis \Ltp{\km}{M} \:.
    \end{aligned}
  \end{equation}
  Thus, $\RHom[\Rhat]{\km}{\Ltp[\Rhat]{D}{\LLam[\fm]{M}}}$ belongs to
  $\Db[\Rhat]$ by the assumption on $M$. By \dcmca{18.2.3} the complex
  $\RHom[\Rhat]{\wN}{D}$ belongs to $\Dfba[\Rhat]$, so
  \thmref{amp-LLam-RHom} yields in view of $(\ast)$ the inequalities,
  \begin{multline*}
    \amp{\Ltpp{\km}{M}} + \cmd[\Rhat]{\RHom[\Rhat]{\wN}{D}} \\
    \dge \amp{\LLam[\fm]{\Ltpp{N}{M}}} \dge
    \cmd[\Rhat]{\RHom[\Rhat]{\wN}{D}} \:.
  \end{multline*}
  Finally, from \dcmca{18.2.31} and faithful flatness of $\Rhat$ as an
  $R$-module one gets:
  \begin{equation*}
    \cmd[\Rhat]{\RHom[\Rhat]{\wN}{D}}  \deq \amp{\wN} \deq \amp{N} \:. \qedhere
  \end{equation*}
\end{prf*}

Nested in the next result, which has a counterpart in
\corref{amp-RGam-Ltp-12}, is Iversen's Amplitude Inequality
\thmcite[3.2]{BIv77}.

\begin{cor}
  \label{cor:amp-LLam-Ltp-12}
  Let $M$ and $N$ be complexes with $M \in \Db$ and $N \in
  \Df$. Assume that $M$ is of finite flat dimension and not
  acyclic. If $M$ is derived $\fm$-complete or belongs to $\Dfb$, then
  there are inequalities,
  \begin{equation*}
    \fd{M} - \inf{M} + \amp{N} \dge \amp{\Ltpp{N}{M}} \dge \amp{N} \:.
  \end{equation*}
\end{cor}

\begin{prf*}
  The complex $M$ belongs to $\Dm$, see \ref{Dm}. If $M$ is derived
  $\fm$-complete, then \prpref{LL}(f) and \eqref{wdt1} yield
  \begin{equation}
    \tag{$\ast$}
    \sup{\Ltpp{\km}{M}} \deq \fd{M} \qqand
    \inf{\Ltpp{\km}{M}} \deq \inf{M} \:;
  \end{equation}
  in particular, the complex $\Ltp{\km}{M}$ belongs to $\Db$. It
  follows from \prpref{LLam} that the complex $\Ltp{N}{M}$ is derived
  $\fm$-complete, so the asserted inequalities follow from
  \thmref{amp-LLam-Ltp} and the equalities $(\ast)$.

  If $M$ belongs to $\Dfb$, then the equalities $(\ast)$ hold by
  \dcmca{16.4.1}, \eqref{wdt1}, and \eqref{Df}. The complex
  $\Ltp{N}{M}$ is in $\Df$, see \dcmca{15.4.3}, so by \eqref{Df} and
  faithful flatness of $\Rhat$ as an $R$-module there are equalities,
  \begin{equation*}
    \amp{\LLam[\fm]{\Ltpp{N}{M}}} \deq \amp{\tpp{\Rhat}{\Ltpp{N}{M}}}
    \deq \amp{\Ltpp{N}{M}} \:.
  \end{equation*}
  Now the asserted inequalities follow from \thmref{amp-LLam-Ltp} and
  the equalities $(\ast)$.
\end{prf*}

Iversen's Amplitude Inequality yields a short proof of Auslander's
zero-divisor conjecture, which is the finitely generated case of the
next example. The other case shows that derived $\fm$-complete modules
of finite flat dimension are ``Auslander modules" in the sense of
Nasehpour \cite{PNs18}.

\begin{exa}
  Let $M$ be an $R$-module and $\pmb{x} = x_1,\ldots,x_n$ an
  $M$-regular sequence in $\fm$. If $M$ has finite flat dimension and
  is finitely generated or derived $\fm$-complete, then $\pmb{x}$ is
  $R$-regular. Indeed, with $K$ denoting the Koszul complex on
  $\pmb{x}$ one has $\amp{\Ltpp{M}{K}} \ge \amp{K}$, so the conclusion
  follows from \dcmca{16.2.31}.
\end{exa}

The inequalities in the next result also hold under slightly different
conditions on the complexes, see \corref{amp-RGam-Ltp-3}.

\begin{cor}
  \label{cor:amp-LLam-Ltp-3}
  Let $M$ and $N$ be complexes with $M \in \Dm$ and $N \in \Dl$. If
  the complex $\Ltp{\km}{M}$ belongs to $\Db$, then there are
  inequalities,
  \begin{equation*}
    \amp{\Ltpp{\km}{M}} + \amp{N} 
    \dge \amp{\Ltpp{N}{M}} \\
    \dge \amp{N} \:.
  \end{equation*}
\end{cor}

\begin{prf*}
  Recall from \eqref{Dlmtor} that the complex $N$ is derived
  $\fm$-torsion, so by \eqref{GLc} and \prpref{LLam} there are
  isomorphisms
  \begin{equation*}
    \Ltp{N}{M} \dqis \Ltp{N}{\LLam[\fm]{M}} \dqis \LLam[\fm]{\Ltpp{N}{M}} \:.
  \end{equation*}
  Now the asserted inequalities follow from \thmref{amp-LLam-Ltp}.
\end{prf*}


\section{Nonvanishing of local cohomology}

\noindent
The first result of this section strengthens \dcmca{18.5.2}, see
\exaref{inf-RGam}. It could be proved in a manner dual to
\thmref{amp-LLam-RHom}, but we obtain it instead by reducing to a
situation where results from \dcmca{18.5} apply.

\begin{thm}
  \label{thm:amp-RGam-Ltp}
  Let $M$ and $N$ be complexes with $M \in \Dm$ and $N\in\Dfba$.  If
  the complex $\Ltp{\km}{M}$ belongs to $\Dba$, then there are
  inequalities,
  \begin{equation*}
    \amp{\Ltpp{\km}{M}} + \cmd{N} 
    \dge
    \amp{\RGam[\fm]{\Ltpp{N}{M}}}
    \dge \cmd{N} \:.
  \end{equation*}
\end{thm}

\begin{prf*}
  Both inequalities are trivial if $N$ is acyclic, so assume that is
  not the case. That is, assume that $N$ belongs to $\Dm$, see
  \ref{Dm}. As $N$ is in $\Dba$ also $\RHom{\km}{N}$ belongs to
  $\Dba$, see \dcmca{7.6.7}, so \eqref{dpt1} and \prpref{ABE}~yield
  \begin{align*}
    \sup{\RGam[\fm]{\Ltpp{N}{M}}} 
    & \deq -\dpt{\Ltpp{N}{M}} \\
    & \deq -\dpt{N}-\dpt{M}+\dptR \:.
  \end{align*}
  Thus, it follows from \corref{ABE} and \eqref{cmd} that it is
  sufficient to establish the inequalities
  \begin{equation}
    \tag{$\ast$}
    \begin{aligned}
      \dim{N} - \wdt{M} &\dge
      -\inf{\RGam[\fm]{\Ltpp{N}{M}}} \\
      &\dge \dim N + \dpt{M} - \dptR \:.
    \end{aligned}
  \end{equation}
  By \eqref{RG} and \eqref{GLc} there are isomorphisms,
  \begin{equation*}
    \RGam[\fm]{\Ltpp{N}{M}} \dqis \Ltp{\RGam[\fm]{N}}{M}
    \dqis \Ltp{\RGam[\fm]{N}}{\LLam[\fm]{M}} \:.
  \end{equation*}
  The left-hand inequality in ($\ast$) now follows from \eqref{GV},
  \eqref{wdt1}, and \dcmca{7.6.8}:
  \begin{equation*}
    \dim{N} - \wdt{M} 
    \deq -\inf{\RGam[\fm]{N}} - \inf{\LLam[\fm]{M}} \dge -\inf{\RGam[\fm]{\Ltpp{N}{M}}} \:.
  \end{equation*}
  To prove the right-hand inequality, first note that one has
  \begin{align*}
    -\inf\RGam[\fm]{\Ltpp{N}{M}} 
    & \dge-\wdt{\RGam[\fm]{\Ltpp{N}{M}}} \\
    & \deq-\wdt{\Ltpp{N}{M}} \\
    & \deq-\wdt{N} - \wdt{M} \\
    & \dge-\wdt{N} - \inf{\RGam[\fm]{M}} - \dimR \:,
  \end{align*}
  where the first inequality holds by \eqref{wdt1} and the second by
  \eqref{wdt3}; the equalities hold by \eqref{wdt2} and
  \eqref{wdt-Ltp}. It follows that the right-hand inequality in
  ($\ast$) holds if one has $\inf{\RGam[\fm]{M}} = -\infty$, so we can
  assume that $\RGam[\fm]{M}$ belongs to $\Db$ as
  $\sup\RGam[\fm]{M}=-\dpt M < \infty$ holds by \eqref{dpt1} and
  \prpref{ABE}. By \prpref{RG}(d) the complex $\RGam[\fm]{M}$ has
  finite flat dimension, so one can apply \dcmca{18.5.1} to get the
  inequality in the computation below. The equalities hold by
  \eqref{RGdim}, \eqref{RG}, and \eqref{dpt2}.
  \begin{align*}
    -\inf{\RGam[\fm]{\Ltpp{N}{M}}}
    & \deq \dim{\RGam[\fm]{\Ltpp{N}{M}}} \\
    & \deq \dim{\Ltpp{N}{\RGam[\fm]{M}}} \\
    & \dge \dim N + \dpt{\RGam[\fm]{M}} - \dptR \\
    & \deq \dim N + \dpt{M} - \dptR \:. \qedhere
  \end{align*}
\end{prf*}

A consequence of Peskine and Szpiro's Intersection Theorem is that
finite length modules of finite projective dimension only exist over
Cohen--Macaulay rings; see \cite[Proposition~6.2.4]{rob}. The next
corollary soups up this classic result, cf.~\eqref{Dlmtor}. The
statement is parallel to \corref{amp-LLam}, but it is worth recalling
that a module is derived $\fm$-torsion if and only if it is
$\fm$-torsion, see \dcmca{13.4.9}.

\begin{cor}
  \label{cor:amp-RGam}
  If there exists a nonzero derived $\fm$-torsion $R$-module of finite
  flat dimension, then $R$ is Cohen--Macaulay.
\end{cor}

\begin{prf*}
  If $M \ne 0$ a derived $\fm$-torsion $R$-module of finite flat
  dimension, then $M$ is in $\Dm$, see \eqref{Dmtor}, and
  $\Ltp{\km}{M}$ is in $\Dba$. As $\amp{\Ltpp{R}{M}} = \amp{M} = 0$
  holds, \thmref{amp-RGam-Ltp} yields $\cmdR = 0$.
\end{prf*}

The next result should be compared with \corref{amp-LLam-Ltp-12}.

\begin{cor}
  \label{cor:amp-RGam-Ltp-12}
  Let $M$ and $N$ be $R$-complexes with $N\in\Dfba$. Assume that $M$
  is of finite flat dimension and not acyclic. If $M$ is derived
  $\fm$-torsion or belongs to $\Dfb$, then there are inequalities,
  \begin{equation*}
    \fd{M} - \wdt{M} + \cmd{N} \dge \cmd{\Ltpp{N}{M}} \dge \cmd{N} \:.
  \end{equation*}
\end{cor}

\begin{prf*}
  The complex $M$ belongs to $\Dm$, see \ref{Dm}. If $M$ is derived
  $\fm$-torsion, then \prpref{RG}(d) and \eqref{wdt1} yield
  \begin{equation}
    \tag{$\ast$}
    \sup{\Ltpp{\km}{M}} \deq  \fd{M} \qqand
    \inf{\Ltpp{\km}{M}} \deq  \wdt{M} \:;
  \end{equation}
  in particular, $\Ltp{\km}{M}$ belongs to $\Dba$. Further, the
  complex $\Ltp{N}{M}$ is derived $\fm$-torsion, see \eqref{LL}, so
  the asserted inequalities follow from \thmref{amp-RGam-Ltp},
  \eqref{Dl}, and the equalities $(\ast)$.

  If $M$ belongs to $\Dfb$, then the equalities $(\ast)$ hold by
  \dcmca{16.4.1} and \eqref{wdt1}. The complex $\Ltp{N}{M}$ belongs to
  $\Dfba$, see \dcmca{15.4.3}, so the asserted inequalities follow
  from \thmref{amp-RGam-Ltp}, \eqref{cmdRG}, and the equalities
  $(\ast)$.
\end{prf*}

Notice that the conclusion in the next result is identical to the
conclusion in \corref{amp-LLam-Ltp-3} but the boundedness conditions
on the complexes are different.

\begin{cor}
  \label{cor:amp-RGam-Ltp-3}
  Let $M$ and $N$ be $R$-complexes with $M \in \Dm$ and $N \in \Dlba$.
  If the complex $\Ltp{\km}{M}$ belongs to $\Dba$, then there are
  inequalities,
  \begin{equation*}
    \amp{\Ltpp{\km}{M}} + \amp{N} \dge \amp{\Ltpp{N}{M}} \dge \amp{N} \:.
  \end{equation*}
\end{cor}

\begin{prf*}
  Recall from \eqref{Dlmtor} that the complex $N$ is derived
  $\fm$-torsion, so $\Ltp{N}{M}$ is per \eqref{RG} derived
  $\fm$-torsion. Now apply \thmref{amp-RGam-Ltp} and~\eqref{Dl}.
\end{prf*}

\begin{exa}
  \label{exa:inf-RGam}
  Let $R$ be Gorenstein of positive Krull dimension and $\fp$ be a
  minimal prime ideal in $R$. Consider the injective $R$-module
  $M = \E{R/\fp} \oplus \E{\km}$. As $R$ is Gorenstein, one has
  $\sup{\Ltpp{\km}{M}} \le \fd{M} < \infty$, so \thmref{amp-RGam-Ltp}
  and \dcmca{18.5.2} both apply to the complex $\Ltp{R}{M} \qis M$.

  One has $\RGam[\fm]{M} = \Gamma_\fm M = \E{\km}$, so
  $\amp{\RGam[\fm]{M}} = 0 = \dpt{M}$ holds. Since the module
  $\E{R/\fp}_\fq \is \E[R_\fq]{R_\fq/\fp_\fq}$ is nonzero for every
  prime ideal $\fq$ in $R$ that contains $\fp$, one has
  $\dim{M} = \dimR$. Therefore, \dcmca{18.5.2} yields
  \begin{equation*}
    \dimR \deq \cmd{\Ltpp{R}{M}} \dge \cmdR \deq  0\:,
  \end{equation*}
  and by the assumption on $R$ the inequality is strict. Yet, the
  right-hand inequality in \thmref{amp-RGam-Ltp} reads $0 \ge 0$ and,
  in fact, so does the left-hand inequality as one has
  \begin{equation*}
    \dimR \dge \fd{M} \dge \sup{\Ltpp{\km}{M}} \dge
    \wdt{M} \deq \dptR \deq \dimR \:,
  \end{equation*}
  by standard (in)equalities, \eqref{wdt1}, and \dcmca{17.4.11 and
    16.2.29}.
\end{exa}

The next result improves \dcmca{13.3.20(c)} for the ideal
$\mathfrak{a} = \fm$.

\begin{prp}
  \label{prp:RGam}
  Let $M$ and $K$ be $R$-complexes with $K \in \Df$. If the complex
  $\RHom{\km}{M}$ belongs to $\Db$, then there is an isomorphism,
  \begin{equation*}
    \RGam[\fm]{\RHom{K}{M}} \dqis \RHom{K}{\RGam[\fm]{M}} \:.
  \end{equation*}
\end{prp}

\begin{prf*}
  The first two isomorphisms in the display below follow from
  \dcmca{13.4.1(d)} and \eqref{LL}. By \prpref{LL}(e) the $R$-complex
  $\LLam[\fm]{M}$ has bounded homology and finite injective dimension.
  As $K$ belongs to $\Df$, the third isomorphism holds by
  \dcmca{13.3.20(c) and 13.4.1(d)}.
  \begin{align*}
    \RGam[\fm]{\RHom{K}{M}} 
    & \dqis \RGam[\fm]{\LLam[\fm]{\RHom{K}{M}}} \\
    & \dqis \RGam[\fm]{\RHom{K}{\LLam[\fm]{M}}} \\
    & \dqis \RHom{K}{\RGam[\fm]{M}} \:. \qedhere
  \end{align*}
\end{prf*}

We derive the next result from \thmref{amp-RGam-Ltp} by way of
Grothendieck Duality.


\begin{thm}
  \label{thm:amp-RGam-RHom}
  Let $M$ and $K$ be complexes with $M \in \Dm$ and $K \in \Df$. If
  the complex $\RHom{\km}{M}$ belongs to $\Db$, then there are
  inequalities,
  \begin{equation*}
    \amp{\RHom{\km}{M}} + \amp{K} \dge \amp{\RGam[\fm]{\RHom{K}{M}}}
    \dge \amp{K} \:.
  \end{equation*}
\end{thm}

\begin{prf*}
  Both inequalities are trivial if $K$ is acyclic, so assume that is
  not the case. That is, assume that $K$ belongs to $\Dm$, see
  \ref{Dm}.  The equalities below hold by \eqref{dpt1} and
  \eqref{dpt-RHom}.
  \begin{equation*}
    \sup{\RGam[\fm]{\RHom{K}{M}}}  
    \deq -\dpt{\RHom{K}{M}} 
    \deq -\wdt{K} - \dpt{M} \:.                                               
  \end{equation*}
  Per \eqref{Df} one has $\wdt{K} = \inf{K}$, so the equalities show
  that both of the asserted inequalities are trivial if
  $\inf{K} = -\infty$ holds.
  
  One can now assume that $K$ belongs to $\Dfbb$. Set
  $\wK = \tp{\Rhat}{K}$; it is by \dcmca{18.3.2(c)} a complex in
  $\Dfbb[\Rhat]$.  The first isomorphism in the display below holds by
  \prpref{RG}(a), and the second follows from \prpref{RGam}. As
  $\RGam[\fm]{M}$ by another application of \prpref{RG}(a) is a
  complex in $\Dmhat$, the third isomorphism is standard, see
  \dcmca{12.3.30}. The ring $\Rhat$ has a dualizing complex, $D$, and
  the fourth isomorphism is Grothendieck Duality, see
  \dcmca{18.2.3}. The last isomorphism holds by \dcmca{12.3.26(c)}, as
  the $\Rhat$-complex $\RGam[\fm]{M}$ has bounded homology and finite
  injective dimension by \prpref{RG}(f,g).
  \begin{align*}
    \RGam[\fm]{\RHom{K}{M}} 
    & \dqis \RGam[\wfm]{\RGam[\fm]{\RHom{K}{M}}} \\
    & \dqis \RGam[\wfm]{\RHom{K}{\RGam[\fm]{M}}} \\
    & \dqis \RGam[\wfm]{\RHom[\Rhat]{\wK}{\RGam[\fm]{M}}} \\    
    & \dqis \RGam[\wfm]{\RHom[\Rhat]{\RHom[\Rhat]{\RHom{\wK}{D}}{D}}{\RGam[\fm]{M}}} \\
    & \dqis \RGam[\wfm]{\Ltpp[\Rhat]{\RHom[\Rhat]{\wK}{D}}{\RHom[\Rhat]{D}{\RGam[\fm]{M}}}} \:.
  \end{align*}
  The complex $\RHom[\Rhat]{D}{\RGam[\fm]{M}}$ belongs to $\Dmhat$,
  see \eqref{Dm}.  The first isomorphism below holds by
  \dcmca{12.3.26(c)}. The second isomorphism holds as one can assume
  that $D$ is normalized, see \dcmca{18.2.23 and 18.2.24}, and the
  third comes from \prpref{RG}(c).
  \begin{equation}
    \tag{$\ast$}
    \begin{aligned}
      \Ltp[\Rhat]{\km}{\RHom[\Rhat]{D}{\RGam[\fm]{M}}} 
      & \dqis \RHom[\Rhat]{\RHom[\Rhat]{\km}{D}}{\RGam[\fm]{M}} \\
      & \dqis \RHom[\Rhat]{\km}{\RGam[\fm]{M}} \\
      & \dqis \RHom{\km}{M} \:.
    \end{aligned}
  \end{equation}
  Thus, $\Ltp[\Rhat]{\km}{\RHom[\Rhat]{D}{\RGam[\fm]{M}}}$ belongs to
  $\Db[\Rhat]$ by the assumption on $M$. By \dcmca{18.2.3} the complex
  $\RHom[\Rhat]{\wK}{D}$ belongs to $\Dfba[\Rhat]$, so
  \thmref{amp-RGam-Ltp} yields in view of $(\ast)$ the inequalities,
  \begin{multline*}
    \amp{\RHom{\km}{M}} + \cmd[\Rhat]{\RHom[\Rhat]{\wK}{D}} \\
    \dge \amp{\RGam[\fm]{\RHom{K}{M}}} \dge
    \cmd[\Rhat]{\RHom[\Rhat]{\wK}{D}} \:.
  \end{multline*}
  Finally, from \dcmca{18.2.31} and faithful flatness of $\Rhat$ as an
  $R$-module one gets
  \begin{equation*}
    \cmd[\Rhat]{\RHom[\Rhat]{\wK}{D}} \deq \amp{\wK} \deq \amp{K} \:. \qedhere
  \end{equation*}
\end{prf*}

A consequence of Grothendieck's Local Duality Theorem is the formula
\begin{equation*}
  \cmd{\RHom{K}{D}} \deq \amp{K}\:,
\end{equation*}
used in the last line of the proof above. The inequalities in the next
corollary can be seen as a generalization, as they collapse to the
formula above for $M=D$ a dualizing $R$-complex. The corollary should
be compared with \corref{amp-LLam-RHom-12}.

\begin{cor}
  \label{cor:amp-RGam-RHom-12}
  Let $M$ and $K$ be complexes with $M \in \Db$ and $K \in
  \Df$. Assume that $M$ is of finite injective dimension and not
  acyclic. If $M$ is derived $\fm$-torsion or belongs to $\Dfb$, then
  there are inequalities,
  \begin{equation*}
    \id{M} - \dpt{M} + \amp{K} \dge \cmd{\RHom{K}{M}}
    \dge \amp{K} \:.
  \end{equation*}
\end{cor}

\begin{prf*}
  The complex $M$ belongs to $\Dm$, see \ref{Dm}. If $M$ is derived
  $\fm$-torsion, then \eqref{dpt1} and \prpref{RG}(f) yield
  \begin{equation}
    \tag{$\ast$}
    \sup{\RHom{\km}{M}} \deq -\dpt{M} \qand
    \inf{\RHom{\km}{M}} \deq -\id{M} \:;
  \end{equation}
  thus the complex $\RHom{\km}{M}$ belongs to $\Db$ as
  $-\dpt{M} \le \sup{M}$ holds by \eqref{dpt1}. It follows from
  \prpref{RGam} that $\RHom{K}{M}$ is derived $\fm$-torsion, so the
  asserted inequalities follow from \thmref{amp-RGam-RHom},
  \eqref{Dl}, and the equalities $(\ast)$.

  If $M$ belongs to $\Dfb$, then the equalities $(\ast)$ hold by
  \eqref{dpt1} and \dcmca{16.4.8}. The complex $\RHom{K}{M}$ belongs
  to $\Df$, see \dcmca{15.4.9}.  Now the asserted inequalities follow
  from \thmref{amp-RGam-RHom}, \eqref{cmdRG}, and the equalities
  $(\ast)$.
\end{prf*}

Note that the conclusion in the next corollary matches the conclusion
in \corref{amp-LLam-RHom-3} but the boundedness conditions on the
complexes are different.

\begin{cor}
  \label{cor:amp-RGam-RHom-3}
  Let $M$ and $K$ be complexes with $M \in \Dm$ and $K \in \Dl$. If
  the complex $\RHom{\km}{M}$ belongs to $\Db$, then there are
  inequalities,
  \begin{equation*}
    \amp{\RHom{\km}{M}} + \amp{K} \dge \amp{\RHom{K}{M}} \dge \amp{K} \:.
  \end{equation*}
\end{cor}

\begin{prf*}
  Recall from \eqref{Dlmtor} that the complex $K$ is derived
  $\fm$-torsion, so by \eqref{GLb} and \prpref{RGam} there are
  isomorphisms
  \begin{equation*}
    \RHom{K}{M} \dqis \RHom{K}{\RGam[\fm]{M}} \dqis \RGam[\fm]{\RHom{K}{M}} \:.
  \end{equation*}
  Now the asserted inequalities follow from \thmref{amp-RGam-RHom}.
\end{prf*}

\begin{exa}
  Let $D$ be a dualizing $R$-complex and $K$ a complex in
  $\Dl$. \corref{amp-RGam-RHom-3} yields
  \begin{equation*}
    \amp{\RHom{K}{D}} \deq \amp{K}\:.
  \end{equation*}
  This equality is well-known, for example from \dcmca{18.2.40 and
    16.1.25(6)}.
\end{exa}

\section{Additional nonvanishing of local (co)homology}

\noindent
In a derived Hom complex one can place conditions of vanishing of
(co)homology with coefficients in $\km$ on either one of the two
arguments. In \thmref{amp-LLam-RHom} the assumption is that the second
complex is of finite $\km$-injective dimension. In the next result,
the assumption is instead that the first complex is of finite
$\km$-flat dimension. We have placed the result in this short closing
section because we obtain it from \thmref{amp-RGam-Ltp} by way of
Matlis Duality.

\begin{thm}
  \label{thm:amp-LLam-RHom-dual}
  Let $M$ and $K$ be complexes with $M \in \Dm$ and $K \in
  \Dartbb$. If the complex $\Ltp{\km}{M}$ belongs to $\Dba$, then
  there are inequalities,
  \begin{equation*}
    \amp{\Ltpp{\km}{M}} + \amp{\LLam[\fm]{K}} 
    \dge \amp{\LLam[\fm]{\RHom{M}{K}}} 
    \dge \amp{\LLam[\fm]{K}} \:.
  \end{equation*}
\end{thm}

\begin{prf*}
  Assume first that $R$ is complete. By Matlis duality \dcmca{18.1.9}
  the complex $N=\Hom{K}{\E{\km}}$ belongs to $\Dfba$ and one has
  \begin{equation}
    \tag{$\ast$}
    K \dqis \RHom{N}{\E{\km}} \:. 
  \end{equation}  
  Thus, by $(\ast)$, Hom--tensor adjunction, \eqref{LL}, and
  \dcmca{16.1.25(6)} there are equalities
  \begin{align*}
    \amp{\LLam[\fm]{\RHom{M}{K}}}
    & \deq\amp{\LLam[\fm]{\RHom{M}{\RHom{N}{\E{\km}}}}} \\
    & \deq \amp{\LLam[\fm]{\RHom{\Ltp{N}{M}}{\E{\km}}}} \\
    & \deq \amp{\RHom{\RGam[\fm]{\Ltpp{N}{M}}}{\E{\km}}} \\
    & \deq \amp{\RGam[\fm]{\Ltpp{N}{M}}} \:.
  \end{align*}
  The asserted inequalities now follow from \thmref{amp-RGam-Ltp} as
  one has
  \begin{align*}
    \cmd N \deq \amp{\RGam[\fm]{N}}
    \deq \amp{\RHom{\RGam[\fm]{N}}{\E{\km}}}
    \deq \amp{\LLam[\fm]{K}}
  \end{align*}
  by \eqref{cmdRG}, \eqref{LL}, and $(\ast)$.
  
  Now consider the general case. Per \prpref{RG}(b,c) and the
  assumptions on $M$, the complex $\RGam[\fm]{M}$ belongs to $\Dmhat$
  and $\Ltp[\Rhat]{\km}{\RGam[\fm]{M}}$ is in $\Dba[\Rhat]$. Further,
  as $K \in \Dart$ one has $\RGam[\fm]{K} \in \Dart[\Rhat]$, i.e.~each
  $\Rhat$-module $\H[n]{\RGam[\fm]{K}}$ is artinian. To see this, it
  suffices to argue that $\H[n]{\RGam[\fm]{K}}$ is artinian as an
  $R$-module; however, in $\D$ one has $\RGam[\fm]{K} \qis K$ as
  $\Dart \subseteq \Dmtor$ by \dcmca{16.1.33}, and hence
  $\H[n]{\RGam[\fm]{K}} \cong \H[n]{K}$, which is an artinian
  $R$-module by assumption.
  
  The first isomorphism in the computation below holds by
  \prpref{LL}(a). The remaining isomorphisms hold by \eqref{LL},
  \eqref{GLb}, and \dcmca{13.4.16(b)}.
  \begin{align*}
    \LLam[\fm]{\RHom{M}{K}}
    & \dqis
      \LLam[\wfm]{\LLam[\fm]{\RHom{M}{K}}}
    \\
    & \dqis
      \LLam[\wfm]{\RHom{\RGam[\fm]{M}}{K}}
    \\
    & \dqis
      \LLam[\wfm]{\RHom{\RGam[\fm]{M}}{\RGam[\fm]{K}}}
    \\
    & \dqis
      \LLam[\wfm]{\RHom[\Rhat]{\RGam[\fm]{M}}{\RGam[\fm]{K}}}
  \end{align*}
  As we have established the desired inequalities over complete local
  rings, we get
  \begin{multline*}
    \amp{\Ltpp[\Rhat]{\km}{\RGam[\fm]{M}}} + \amp{\LLam[\wfm]{\RGam[\fm]{K}}} \\
    \dge \amp{\LLam[\fm]{\RHom{M}{K}}} \dge
    \amp{\LLam[\wfm]{\RGam[\fm]{K}}} \:.
  \end{multline*}
  By \prpref{RG}(c) one has
  $\amp{\Ltpp[\Rhat]{\km}{\RGam[\fm]{M}}} = \amp{\Ltpp{\km}{M}}$ and,
  finally,
  \begin{equation*}
    \LLam[\wfm]{\RGam[\fm]{K}}
    \dqis
    \LLam[\fm]{\RGam[\fm]{K}}
    \dqis
    \LLam[\fm]{K}
  \end{equation*}
  by independence of base for local homology, see \dcmca{13.1.21(a),
    and 13.4.1(c)}.
\end{prf*}

\begin{cor}
  \label{cor:amp-LLam-RHom-dual-12}
  Let $M$ and $K$ be $R$-complexes with $K \in \Dartbb$. If $M$ is
  derived $\fm$-torsion of finite flat dimension and not acyclic, then
  there are inequalities,
  \begin{equation*}
    \fd{M} - \wdt{M} + \amp{\LLam[\fm]{K}} 
    \dge \amp{\RHom{M}{K}} 
    \dge \amp{\LLam[\fm]{K}} \:.
  \end{equation*}  
\end{cor}

\begin{prf*}
  The complex $M$ belongs to $\Dm$, see \ref{Dm}, and \prpref{RG}(d)
  and \eqref{wdt1} yield
  \begin{equation*}
    \sup{\Ltpp{\km}{M}} \deq \fd{M} \qqand \inf{\Ltpp{\km}{M}} \deq \wdt{M} \:;
  \end{equation*}
  in particular, the complex $\Ltp{\km}{M}$ belongs to
  $\Dba$. Further, per \eqref{LL} the complex $\RHom{M}{K}$ is derived
  $\fm$-complete, so the asserted inequalities follow from
  \thmref{amp-LLam-RHom-dual} and the equalities displayed above.
\end{prf*}

\begin{cor}
  \label{cor:amp-LLam-RHom-dual-3}
  Let $M$ and $K$ be complexes with $M \in \Dm$ and $K \in \Dlbb$. If
  the complex $\Ltp{\km}{M}$ belongs to $\Dba$, then there are
  inequalities,
  \begin{equation*}
    \amp{\Ltpp{\km}{M}} + \amp{K} 
    \dge \amp{\RHom{M}{K}}
    \dge \amp{K} \:.
  \end{equation*}  
\end{cor}

\begin{prf*}
  The complex $K$ is derived $\fm$-complete by \eqref{Dlmtor} and
  belongs to $\Dartbb$. In particular, the complex $\RHom{M}{K}$ is
  derived $\fm$-complete by \eqref{LL}, so the asserted inequalities
  follow from \thmref{amp-LLam-RHom-dual}.
\end{prf*}

The next result is dual to \thmref{amp-RGam-RHom}; it applies to show
that the inequalities in \corref{amp-LLam-RHom-dual-3} also hold under
the slightly different assumptions $M \in \Dm$ with
$\Ltp{\km}{M} \in \Db$ and $K \in \Dl$. This is the phenomenon
observed in \corref[Corollaries~]{amp-LLam-RHom-3} and
\corref[]{amp-RGam-RHom-3} as well as in
\corref[Corollaries~]{amp-LLam-Ltp-3} and \corref[]{amp-RGam-Ltp-3}.
  
\begin{thm}
  \label{thm:amp-RGam-RHom-dual}
  Let $M$ and $K$ be complexes with $M \in \Dm$ and $K \in \Dart$. If
  the complex $\Ltp{\km}{M}$ belongs to $\Db$, then there are
  inequalities,
  \begin{equation*}
    \amp{\Ltpp{\km}{M}} + \amp{K}
    \dge \amp{\RGam[\fm]{\RHom{M}{K}}} 
    \dge \amp{K} \:.
  \end{equation*}
\end{thm}

\begin{prf*}
  Using the same technique as in the proof of
  \thmref{amp-LLam-RHom-dual}, this result is derived from
  \thmref{amp-RGam-RHom} by way of Matlis Duality.
\end{prf*}

\begin{cor}
  \label{cor:amp-RGam-RHom-dual}
  Let $M$ and $K$ be complexes with $M \in \Dfb$ and $K \in \Dart$. If
  $M$ has finite projective dimension and is not acyclic, then there
  are inequalities,
  \begin{equation*}
    \pd{M} - \inf{M} + \amp{K}
    \dge \amp{\RHom{M}{K}} 
    \dge \amp{K} \:.
  \end{equation*}
\end{cor}

\begin{prf*}  
  As $M$ belongs to $\Dfb$, there are by \dcmca{16.4.1} and \eqref{Df}
  equalities,
  \begin{equation*}
    \sup{\Ltpp{\km}{M}} \deq \pd{M} \qqand \inf{\Ltpp{\km}{M}} \deq \inf{M} \:.
  \end{equation*}
  One has
  \mbox{$\supp{\RHom{M}{K}} = \supp{M} \cap \supp{K} \subseteq
    \set{\fm}$} by \dcmca{17.1.11(a) and 16.1.33}, so the complex
  $\RHom{M}{K}$ is derived $\fm$-torsion, see \eqref{Dmtor}. Now the
  asserted inequalities hold by \thmref{amp-RGam-RHom-dual} and the
  equalities displayed above.
\end{prf*}


\begin{thebibliography}{10}

\bibitem{AJL-97} Leovigildo Alonso~Tarr\'io, Ana Jerem\'ias~L\'opez,
  and Joseph Lipman, \emph{Local homology and cohomology on schemes},
  Ann. Sci. \'Ecole Norm. Sup.  (4) \textbf{30} (1997), no.~1,
  1--39. \MR{MR1422312}

\bibitem{YAn18} Yves Andr\'{e}, \emph{Perfectoid spaces and the
    homological conjectures}, Proceedings of the {I}nternational
  {C}ongress of {M}athematicians---{R}io de {J}aneiro
  2018. {V}ol. {II}. {I}nvited lectures, World Sci. Publ., Hackensack,
  NJ, 2018, pp.~277--289. \MR{MR3966766}

\bibitem{HBs63} Hyman Bass, \emph{On the ubiquity of {G}orenstein
    rings}, Math. Z. \textbf{82} (1963), 8--28. \MR{MR0153708}

\bibitem{BIK-12} David~J. Benson, Srikanth~B. Iyengar, and Henning
  Krause, \emph{Colocalizing subcategories and cosupport}, J. Reine
  Angew. Math. \textbf{673} (2012), 161--207. \MR{MR2999131}

\bibitem{dcmca} Lars~Winther Christensen, Hans-Bj{\o}rn Foxby, and
  Henrik Holm, \emph{Derived {C}ategory {M}ethods in {C}ommutative
    {A}lgebra}, Springer Monographs in Mathematics, Springer, Cham,
  2024. \MR{MR4890472}

\bibitem{HBF79} Hans-Bj{\o}rn Foxby, \emph{Bounded complexes of flat
    modules}, J. Pure Appl.  Algebra \textbf{15} (1979), no.~2,
  149--172. \MR{MR0535182}

\bibitem{hha} Hans-Bj{\o}rn Foxby, \emph{Hyperhomological algebra \&
    commutative rings}, lecture notes, 1998.

\bibitem{HBFSIn03} Hans-Bj{\o}rn Foxby and Srikanth Iyengar,
  \emph{Depth and amplitude for unbounded complexes}, Commutative
  algebra ({G}renoble/{L}yon, 2001), Contemp.  Math., vol. 331,
  Amer. Math. Soc., Providence, RI, 2003, pp.~119--137.
  \MR{MR2013162}

\bibitem{AFrPJr03b} Anders Frankild and Peter J{\o}rgensen,
  \emph{Homological identities for differential graded algebras},
  J.~Algebra \textbf{265} (2003), no.~1, 114--135. \MR{MR1984902}

\bibitem{JPGJPM92} J.~P.~C. Greenlees and J.~P. May, \emph{Derived
    functors of {$I$}-adic completion and local homology}, J. Algebra
  \textbf{149} (1992), no.~2, 438--453. \MR{MR1172439}

\bibitem{lch} Robin Hartshorne, \emph{Local cohomology}, A seminar
  given by A.~Grothendieck, Harvard University, Fall 1961, vol.~41,
  Springer-Verlag, Berlin, 1967.  \MR{MR0224620}

\bibitem{BIv77} Birger Iversen, \emph{Amplitude inequalities for
    complexes}, Ann. Sci. \'Ecole Norm. Sup. (4) \textbf{10} (1977),
  no.~4, 547--558. \MR{MR0568903}

\bibitem{PNs18} Peyman Nasehpour, \emph{Auslander modules},
  Beitr. Algebra Geom. \textbf{59} (2018), no.~4,
  617--624. \MR{MR3871097}

\bibitem{CPsLSz73} C.~Peskine and L.~Szpiro, \emph{Dimension
    projective finie et cohomologie locale. {A}pplications \`a la
    d\'{e}monstration de conjectures de {M}.  {A}uslander, {H}. {B}ass
    et {A}. {G}rothendieck}, Inst. Hautes \'{E}tudes
  Sci. Publ. Math. (1973), no.~42, 47--119. \MR{MR0374130}

\bibitem{PSY-14} Marco Porta, Liran Shaul, and Amnon Yekutieli,
  \emph{On the homology of completion and torsion},
  Algebr. Represent. Theory \textbf{17} (2014), no.~1,
  31--67. \MR{MR3160712}

\bibitem{PRb87} Paul Roberts, \emph{Le th\'eor\`eme d'intersection},
  C. R. Acad. Sci. Paris S\'er. I Math. \textbf{304} (1987), no.~7,
  177--180. \MR{MR0880574}

\bibitem{rob} Paul~C. Roberts, \emph{Multiplicities and {C}hern
    classes in local algebra}, Cambridge Tracts in Mathematics,
  vol. 133, Cambridge University Press, Cambridge,
  1998. \MR{MR1686450}

\bibitem{PSc03} Peter Schenzel, \emph{Proregular sequences, local
    cohomology, and completion}, Math. Scand. \textbf{92} (2003),
  no.~2, 161--180. \MR{MR1973941}

\bibitem{AMS90} Anne-Marie Simon, \emph{Some homological properties of
    complete modules}, Math.  Proc. Cambridge
  Philos. Soc. \textbf{108} (1990), no.~2, 231--246.  \MR{MR1074711}

\bibitem{PScAMS} Peter Schenzel and Anne-Marie Simon,
  \emph{Completion, \v{C}ech and local homology and cohomology},
  Springer Monographs in Mathematics, Springer, Cham,
  2018. \MR{MR3838396}

\bibitem{JRS90} Jan~R. Strooker, \emph{Homological questions in local
    algebra}, London Mathematical Society Lecture Note Series,
  vol. 145, Cambridge University Press, Cambridge,
  1990. \MR{MR1074178}

\bibitem{SYs98b} Siamak Yassemi, \emph{Width of complexes of modules},
  Acta Math. Vietnam.  \textbf{23} (1998), no.~1,
  161--169. \MR{MR1628029}

\end{thebibliography}

\def\soft#1{\leavevmode\setbox0=\hbox{h}\dimen7=\ht0\advance \dimen7
  by-1ex\relax\if t#1\relax\rlap{\raise.6\dimen7
    \hbox{\kern.3ex\char'47}}#1\relax\else\if T#1\relax
  \rlap{\raise.5\dimen7\hbox{\kern1.3ex\char'47}}#1\relax \else\if
  d#1\relax\rlap{\raise.5\dimen7\hbox{\kern.9ex
      \char'47}}#1\relax\else\if D#1\relax\rlap{\raise.5\dimen7
    \hbox{\kern1.4ex\char'47}}#1\relax\else\if l#1\relax
  \rlap{\raise.5\dimen7\hbox{\kern.4ex\char'47}}#1\relax \else\if
  L#1\relax\rlap{\raise.5\dimen7\hbox{\kern.7ex
      \char'47}}#1\relax\else\message{accent \string\soft \space #1
    not defined!}#1\relax\fi\fi\fi\fi\fi\fi}
\providecommand{\MR}[1]{\mbox{\href{http://www.ams.org/mathscinet-getitem?mr=#1}{#1}}}
\renewcommand{\MR}[1]{\mbox{\href{http://www.ams.org/mathscinet-getitem?mr=#1}{#1}}}
\providecommand{\arxiv}[2][AC]{\mbox{\href{http://arxiv.org/abs/#2}{\sf
      arXiv:#2 [math.#1]}}} \def\cprime{$'$}
\providecommand{\bysame}{\leavevmode\hbox to3em{\hrulefill}\thinspace}
\providecommand{\MR}{\relax\ifhmode\unskip\space\fi MR }
\providecommand{\MRhref}[2]{%
  \href{http://www.ams.org/mathscinet-getitem?mr=#1}{#2} }
\providecommand{\href}[2]{#2}

\end{document}